\documentclass[10pt,righttag]{amsart}
\usepackage{amssymb,verbatim,amscd,mathrsfs,stmaryrd,wasysym,latexsym, bbm}
\usepackage[all]{xy} 
\usepackage{cancel,color, multirow}
\usepackage{lscape}
\begin{document}

\newcommand{\V}{{\mathcal V}}     
\renewcommand{\O}{{\mathcal O}}
\newcommand{\LL}{\mathcal L}
\newcommand{\mcO}{\mathcal{O}}
\newcommand{\mcF}{\mathcal{F}}
\newcommand{\mcL}{\mathcal{L}}
\newcommand{\mCG}{\mathcal{G}}
\newcommand{\mcH}{\mathcal{H}}
\newcommand{\mcI}{\mathcal{I}}
\newcommand{\mcK}{\mathcal{K}}

\newcommand{\Ext}{\hbox{\rm Ext}}
\newcommand{\Tor}{\hbox{\rm Tor}}
\newcommand{\Hom}{\hbox{Hom}}
\newcommand{\Spec}{\hbox{Spec }}
\newcommand{\Proj}{\hbox{Proj }}
\newcommand{\Mod}{\hbox{Mod}}
\newcommand{\GrMod}{\hbox{GrMod}}
\newcommand{\grmod}{\hbox{gr-mod}}
\newcommand{\Tors}{\hbox{Tors}}
\newcommand{\gr}{\hbox{gr}}
\newcommand{\tors}{\hbox{tors}}
\newcommand{\rank}{\hbox{rank}}
\newcommand{\End}{\hbox{{\rm End}}}
\newcommand{\Der}{\hbox{Der}}
\newcommand{\GKdim}{\hbox{GKdim}}
\newcommand{\gldim}{\hbox{gldim}}
\newcommand{\im}{\hbox{im}}
\renewcommand{\ker}{\hbox{ker}}
\newcommand{\Char}{\hbox{char}}
\newcommand{\colim}{\hbox{colim}}
\newcommand{\depth}{\hbox{depth}}
\def\bee{\begin{eqnarray}}
\def\eee{\end{eqnarray}}

\newcommand{\AGr}{\hbox{A-Gr}}
\newcommand{\lonto}{{\protect \longrightarrow\!\!\!\!\!\!\!\!\longrightarrow}}

\renewcommand{\c}{\cancel}
\newcommand{\fh}{\frak h}  
\newcommand{\fp}{\frak p}
\newcommand{\fq}{\frak q}
\newcommand{\fr}{\frak r}
\newcommand{\mf}{\mathfrak}
\newcommand{\m}{{\mu}}
\newcommand{\gl}{{\frak g}{\frak l}}
\newcommand{\ssl}{{\frak s}{\frak l}}
\newcommand{\tw}{{\rm tw}}

\newcommand{\ds}{\displaystyle}
\newcommand{\s}{\sigma}
\renewcommand{\l}{\lambda}
\renewcommand{\a}{\alpha}
\renewcommand{\b}{\beta}
\newcommand{\G}{\Gamma}
\newcommand{\g}{\gamma}
\newcommand{\z}{\zeta}
\newcommand{\e}{\epsilon}
\renewcommand{\d}{\delta}
\newcommand{\p}{\rho}
\renewcommand{\t}{\tau}
\newcommand{\n}{\nu}
\newcommand{\x}{\chi}
\newcommand{\w}{\omega}
\renewcommand{\i}{\iota}

\newcommand{\A}{{\Bbb A}}
\newcommand{\C}{{\Bbb C}}
\newcommand{\N}{{\Bbb N}}
\newcommand{\Z}{{\Bbb Z}}
\newcommand{\ZZ}{{\Bbb Z}}
\newcommand{\Q}{{\Bbb Q}}
\renewcommand{\k}{\mathbb K}

\newcommand{\E}{{\mathcal E}}
\newcommand{\K}{{\mathcal K}}
\renewcommand{\L}{{\mathcal L}}
\renewcommand{\S}{{\mathcal S}}
\newcommand{\T}{{\mathcal T}}

\newcommand{\GL}{{GL}}

\newcommand{\rowxy}{(x\ y)}
\newcommand{\colxy}{ \left({\begin{array}{c} x \\ y \end{array}}\right)}
\newcommand{\scolxy}{\left({\begin{smallmatrix} x \\ y
\end{smallmatrix}}\right)}

\renewcommand{\P}{{\Bbb P}}

\newcommand{\la}{\langle}
\newcommand{\ra}{\rangle}
\newcommand{\tensor}{\otimes}
\newcommand{\tsr}{\tensor}
\newcommand{\ol}{\overline}

\newtheorem{thm}{Theorem}[section]
\newtheorem{lemma}[thm]{Lemma}
\newtheorem{cor}[thm]{Corollary}
\newtheorem{prop}[thm]{Proposition}
\newtheorem{claim}[thm]{Claim}

\theoremstyle{definition}
\newtheorem{defn}[thm]{Definition}
\newtheorem{notn}[thm]{Notation}
\newtheorem{ex}[thm]{Example}
\newtheorem{rmk}[thm]{Remark}
\newtheorem{rmks}[thm]{Remarks}
\newtheorem{note}[thm]{Note}
\newtheorem{example}[thm]{Example}
\newtheorem{problem}[thm]{Problem}
\newtheorem{ques}[thm]{Question}
\newtheorem{conj}[thm]{Conjecture}
\newtheorem{thingy}[thm]{}

\newcommand{\onto}{{\protect \rightarrow\!\!\!\!\!\rightarrow}}
\newcommand{\donto}{\put(0,-2){$|$}\put(-1.3,-12){$\downarrow$}{\put(-1.3,-14.5) 

{$\downarrow$}}}

\newcounter{letter}
\renewcommand{\theletter}{\rom{(}\alph{letter}\rom{)}}

\newenvironment{lcase}{\begin{list}{~~~~\theletter} {\usecounter{letter}
\setlength{\labelwidth4ex}{\leftmargin6ex}}}{\end{list}}

\newcounter{rnum}
\renewcommand{\thernum}{\rom{(}\roman{rnum}\rom{)}}

\newenvironment{lnum}{\begin{list}{~~~~\thernum}{\usecounter{rnum}
\setlength{\labelwidth4ex}{\leftmargin6ex}}}{\end{list}}

\thispagestyle{empty}

\title[Quantum $\P^2$s as graded twisted tensor products]{Quantum projective planes as certain graded twisted tensor products}

\keywords{Artin-Schelter regular algebras, noncommutative algebraic geometry, twisted tensor products, Sklyanin algebras}

\author[  Conner, Goetz ]{ }

  \subjclass[2010]{16S37, 16W50}
\maketitle

\begin{center}

\vskip-.2in Andrew Conner \\
\bigskip

Department of Mathematics and Computer Science\\
Saint Mary's College of California\\
Moraga, CA 94575\\
\bigskip

 Peter Goetz \\
\bigskip

Department of Mathematics\\ Humboldt State University\\
Arcata, California  95521
\\ \ \\

\end{center}

\setcounter{page}{1}

\thispagestyle{empty}

\vspace{0.2in}

\begin{abstract}

Let $\k$ be an algebraically closed field. Building upon previous work, we classify, up to isomorphism of graded algebras, quadratic graded twisted tensor products of $\k[x,y]$ and $\k[z]$. When such an algebra is Artin-Schelter regular, we identify its point scheme and type, in the sense of \cite{IM}. We also describe which three-dimensional Sklyanin algebras contain a subalgebra isomorphic to a quantum $\P^1$, and we show that every algebra in this family is a graded twisted tensor product of $\k_{-1}[x,y]$ and $\k[z]$.
\end{abstract}

\bigskip

\section{Introduction}
\label{introduction}

In \cite{Cap}, \v{C}ap, Schichl, and Van\v{z}ura were led to define the notion of a twisted tensor product of algebras over a field as an answer to the following question from noncommutative differential geometry: 

\begin{quotation}
Given two algebras that represent noncommutative spaces, what is an appropriate representative of a noncommutative product of those spaces? 
\end{quotation}

This paper is motivated by a related question in noncommutative projective geometry. If $A$ is a noetherian graded $\k$-algebra, finitely generated in degree 1, the noncommutative projective scheme $(\Proj A, \mcO)$ is, by definition, the quotient category $\Proj A ={\rm GrMod}(A)/{\rm Tors}(A)$ of graded right modules, modulo torsion modules, with $\mcO$ the image of $A_A$ in $\Proj A$. In light of Serre's theorem in commutative algebraic geometry, the category $\Proj A$  can be considered as the category of quasi-coherent sheaves on an imagined noncommutative scheme. 
 In this context, we say $\Proj A$ (or just $A$) is a \emph{quantum projective space} if $A$ is a  quadratic Artin-Schelter regular algebra; more precisely, if $A$ has global dimension $n+1$, then we say that $\Proj A$ (or $A$) is a quantum $\P^n$. We raise the following natural question: 
 
 \begin{quotation}
Which quantum projective spaces are isomorphic to graded twisted tensor products of other quantum projective spaces?
\end{quotation}

It is not hard to show that every quantum $\P^1$ is a graded twisted tensor product of $\k[x]$ and $\k[y]$. One of our main results in this paper is the determination of all of the quantum $\P^2$'s that are isomorphic as graded algebras to graded twisted tensor products of the ``classical'' $\P^1$, $\k[x,y]$, with $\k[z]$. This builds upon results in \cite{CG3}, where all graded twisted tensor products of $\k[x,y]$ with $\k[z]$ were classified up to a notion of equivalence stronger than that of graded algebra isomorphism. 

An intriguing aspect of Artin-Schelter regular algebras is the extent to which their algebraic properties are encoded geometrically. To every quantum $\P^2$, one can associate a subscheme $E\subset \P^2$ and an automorphism $\s:E\to E$ which determine the algebra up to isomorphism. In characterizing which quantum $\P^2$'s are graded twisted tensor products of $\k[x,y]$ and $\k[z]$, we also determine this geometric data. 

The first classification results on quantum $\P^2$ are due to Artin and Schelter in \cite{AS} and Artin-Tate-van den Bergh in \cite{ATVI}. These results provide a ``generic'' classification up to algebra isomorphism which has subsequently been refined through the work of many authors; we refer the reader to \cite{IM} for details. In characteristic 0, a classification of three-dimensional quadratic Artin-Schelter regular algebras up to isomorphism and Morita equivalence has recently been completed \cite{IM, Mat} using the theory of \emph{geometric algebras} developed in \cite{Mori} and \cite{MU}. Though we do not assume ${\rm char}\ \k=0$ in this paper, we adopt the naming conventions introduced by these authors, in which the types of algebras are differentiated according to the geometric data $(E,\s)$, see Section \ref{preliminaries}.  Detailed mappings from the classification of \cite{CG3} to types from \cite{IM} are provided in Sections 3, 4, and 5.

\begin{thm}[Theorem \ref{point scheme of T(g, h)}, Theorem \ref{ReducibleSummary}, Theorem \ref{OreSummary}]
\label{introTypes1}
Assume ${\rm char}\ \k \neq 2$. Let $T$ be an Artin-Schelter regular graded twisted tensor product of $\k[x,y]$ and $\k[z]$. 
\begin{enumerate}
\item If $T$ is of type $EC$, then the point scheme $E$ is an elliptic curve, and the automorphism $\s$ is multiplication by $-1$ in the group law on $E$. ($T$ is of type $B$.)
\item If $T$ is not of type $EC$, then $T$ belongs to one of the following types:
$$P_1, P_2, S_1, S_1', S_2, S_2', T_1, T', WL_1, WL_2, WL_3, TL_1, TL_2, TL_4, NC_2.$$ 
\end{enumerate}
 Moreover, each type listed above is the type of an Artin-Schelter regular graded twisted tensor product of $\k[x,y]$ and $\k[z]$. 
\end{thm}

Comparing Theorem \ref{introTypes1} to the results of \cite{IM, Mat} shows that, in characteristic 0, not all quantum $\P^2$'s can be expressed as graded twisted tensor products of $\k[x,y]$ and $\k[z]$, even up to graded Morita equivalence. We summarize this comparison as follows.

\begin{thm}
\label{introTypes2}
Assume ${\rm char}\ \k =0$. Let $A$ be a three-dimensional quadratic Artin-Schelter regular algebra. 
\begin{enumerate}
\item (Remarks  \ref{T(g,h) remark}, \ref{complete reducible cases}, \ref{complete Ore types}) If $A$ belongs to one of the following types: $$EC \text{ (subtype $B$), } S'_2, T_1, WL_2, WL_3, TL_4, NC_2, $$ then $A$ is isomorphic to a graded twisted tensor product of $\k[x,y]$ and $\k[z]$.

\item (Remark \ref{complete Ore cases}) If $A$ belongs to one of the following types: $$P, S, S', T, T', WL, TL,$$ then $A$ is graded Morita equivalent to a graded twisted tensor product of $\k[x,y]$ and $\k[z]$. 

\item If $A$ is of type $CC$, then $A$ is not graded Morita equivalent to a graded twisted tensor product of $\k[x,y]$ and $\k[z]$. 
\end{enumerate}
\end{thm}

Not every Morita equivalence class of quantum $\P^2$'s contains a graded twisted tensor product of $\k[x,y]$ and $\k[z]$. By \cite[Theorem 3.2]{IM}, the only ones that do not are of Type NC, CC or EC. However, the unique isomorphism class of Type CC is a graded Ore extension of the Jordan plane, and every algebra of Type $NC_1$ is a graded Ore extension of a skew polynomial algebra $\k_{a}[x,y]$. Since graded Ore extensions are twisted tensor products, Theorem \ref{introTypes2} and the preceding observations imply that every quantum $\P^2$ that is not of type EC is graded Morita equivalent to a twisted tensor product (in fact, a graded Ore extension) of \emph{some} quantum $\P^1$ with $\k[z]$.

In the last section of the paper,  we prove there are quantum $\P^2$'s that are not (Morita equivalent to) graded twisted tensor products of \emph{any} quantum $\P^1$ with $\k[z]$. Our results in this direction concern the three-dimensional Sklyanin algebras $S(a,b,c)$, see Section \ref{Sklyanins} for the relevant definitions. We argue first, using noncommutative algebraic geometry, that $\k_{-1}[x,y]$ is the only quantum $\P^1$ that can arise as a subalgebra of $S(a,b,c)$, when the point scheme of $S(a,b,c)$ is an elliptic curve, and this can only occur when $a=b$.

\begin{thm}[Theorem \ref{not a skew polynomial ring}, Proposition \ref{not a Jordan plane}]
\label{introSklyanin1}
Assume ${\rm char}\ \k\neq 2,3$. Let $S(a,b,c)$ be a three-dimensional Sklyanin algebra of type EC.\ Then $S(a,b,c)$ contains a subalgebra isomorphic to a quantum $\P^1$ if and only if $a=b$. The only quantum $\P^1$ contained in $S(a,a,c)$ is $\k_{-1}[x,y]$.
 \end{thm}

Thus, when ${\rm char}\ \k\neq 2,3$, a Sklyanin algebra of type EC is graded Morita equivalent to a graded twisted tensor product of a quantum $\P^1$ with $\k[z]$ only if it is graded Morita equivalent to $S(1,1,c)$ where $c\neq 0$ and $(3c)^3\neq (2+c^3)^3$.

Finally, we show that, in characteristic 0, every three-dimensional Sklyanin algebra of type EC that contains a subalgebra isomorphic to $\k_{-1}[x,y]$ is a graded twisted tensor product of $\k_{-1}[x,y]$ and $\k[z]$.
 
\begin{thm}[Theorem \ref{P(a,b) is an S(1,1,c)}]
\label{introSklyanin2}
Assume ${\rm char}\ \k=0$. A three-dimensional Sklyanin algebra $S(1,1,c)$ of type EC is isomorphic as a graded algebra to a twisted tensor product $P(a)$ for some $a\in \k-\{0,1\}$. 
\end{thm}

The outline of the paper is as follows. Definitions and notation related to twisted tensor products and types of Artin-Schelter regular algebras are given in Section 2. In particular, the three main families of graded twisted tensor products from \cite{CG3} are described. In Sections 3, 4, and 5, algebras from these three families are classified up to isomorphism, and the associated geometric data is computed. Finally, results related to the three-dimensional Sklyanin algebras of type EC are the subject of Section 6.

\section{Preliminaries}
\label{preliminaries}

In this paper we work with algebras over a field $\k$. We assume throughout that $\k$ is algebraically closed and ${\rm char}\, \k\neq 2$. The tensor product of $\k$-vector spaces $V$ and $W$ is denoted $V\tsr W$, and when $V$ and $W$ are $\N$-graded, the space $V\tsr W$ is graded by the formula
$$(V\tsr W)_m=\bigoplus_{k+\ell=m} V_k\tsr W_{\ell}.$$ We write $V^*$ for the linear dual of the $\k$-vector space $V$.

By a \emph{graded $\k$-algebra}, we mean a connected, $\N$-graded, locally finite-dimensional $\k$-algebra, generated in degree 1. We require morphisms of graded algebras and vector spaces to preserve degree. 
We call a graded $\k$-algebra $A$ \emph{quadratic} if there exists a finite dimensional $\k$-vector space $V$ and a subspace $R\subseteq V\tsr V$ such that $A\cong T(V)/\la R\ra$ as graded algebras, where $T(V)$ is the tensor algebra.

\subsection{Graded twisted tensor products}
Let $A$ and $B$ be graded $\k$-algebras. A \emph{graded twisted tensor product} of $A$ and $B$ is a triple $(C,i_A, i_B)$ where $C$ is a graded $\k$-algebra and $i_A:A\to C$ and $i_B:B\to C$ are graded algebra inclusions such that the $\k$-linear map $A\tsr B\to C$ given by $a\tsr b\mapsto i_A(a)i_B(b)$ is an isomorphism. In this paper we are primarily concerned with the setting where $A=\k[x,y]$, $B=\k[z]$, and $C$ is a quadratic algebra.

A graded twisted tensor product imparts the structure of an associative algebra to $A\tsr B$, implicitly determining a product map $(A\tsr B)^{\tsr 2}\to A\tsr B$. The product can be described somewhat more directly via the notion of a twisting map. Given any $\k$-linear map $\t:B\tsr A\to A\tsr B$, one can define $\mu_{\t}:(A\tsr B)^{\tsr 2}\to A\tsr B$ by $\mu_{\t}=(\mu_A\tsr\mu_B)(1\tsr \t\tsr 1)$. However, $\mu_{\t}$ does not define an associative product, in general.

A (unital) \emph{graded twisting map} for $A$ and $B$ is a graded $\k$-linear map $\t:B\tsr A\to A\tsr B$ such that  $\t(1\tsr a)=a\tsr 1$, $\t(b\tsr 1)=1\tsr b$, and 
$$\t(\mu_B\tsr \mu_A)=(\mu_A\tsr \mu_B)(1\tsr \t\tsr 1)(\t\tsr\t)(1\tsr \t\tsr 1)$$
where $\mu_A$ and $\mu_B$ denote the product maps on $A$ and $B$, respectively. This condition is precisely what is required for $\mu_{\t}$ to satisfy associativity.

\begin{prop}{\cite[Proposition 2.3]{C-G}}
Let $A$ and $B$ be graded algebras and $\t:B\tsr A\to A\tsr B$ a graded $\k$-linear map. Then $\t$ is a graded twisting map if and only if $(A\tsr B,\mu_{\t})$ is an associative graded $k$-algebra.
\end{prop}

In particular, if $\t:B\tsr A\to A\tsr B$ is a graded twisting map, the graded algebra $A\tsr_{\t} B=(A\tsr B,\mu_{\t})$, together with the canonical inclusions of $A$ and $B$, is a graded twisted tensor product of $A$ and $B$. If $A\tsr_{\t} B\cong A\tsr_{\t'} B$ as graded $\k$-algebras, it does not necessarily follow that $\t=\t'$. In \cite{Cap}, the authors define a notion of isomorphism  that uniquely identifies twisted tensor products by their twisting maps. In \cite{CG3}, we classified quadratic twisted tensor products $\k[x,y]\tsr_{\t} \k[z]$ up to a notion of equivalence weaker than that of \cite{Cap}, yet stronger than graded algebra isomorphism. The precise definitions of these notions are not needed in this paper.

By \cite[Theorem 1.2]{C-G}, a quadratic twisted tensor product $\k[x,y]\tsr_{\t} \k[z]$ is uniquely determined by $\t(z\tsr x)$ and $\t(z\tsr y)$, and by \cite[Proposition 2.5]{C-G} we have 
$$C=\k[x,y]\tsr_{\t}\k[z]\cong \k\la x,y,z\ra/\la xy-yx, zx-\t(zx), zy-\t(zy)\ra.$$
Here $\t(zx)$ and $\t(zy)$ are the images of $\t(z\tsr x)$ and $\t(z\tsr y)$ in $\k\la x,y,z\ra$ under the obvious identification with $T(C_1)$.  A main result from \cite{CG3} is the following, which does not require any restrictions on the characteristic of $\k$.

\begin{thm}{\cite[Theorem 1.4]{CG3}}
\label{old classification}
Let $\k$ be an algebraically closed field. A quadratic twisted tensor product of $\k[x,y]$ and $\k[z]$ is equivalent to one determined by 
\begin{align*}
\t(zx)&=ax^2+bxy+cy^2+dxz+eyz+fz^2\\
\t(zy)&=Ax^2+Bxy+Cy^2+Dyz
\end{align*}
where $e, f, A\in \{0,1\}$. 
\end{thm}

We alert the reader to a small notational shift from \cite{CG3}: in \cite[Theorem 1.4]{CG3} we wrote $E$ instead of $D$ in the $\t(zy)$ equation. In this paper $E$ always refers to a point scheme, as defined below. When comparing results of this paper with corresponding results in \cite{CG3} one should replace $D$ with $E$.

We divided the algebras of Theorem \ref{old classification} into three families: \emph{Ore type} (when $f=0$), \emph{reducible type} (when $f=1, A=0$) and \emph{elliptic type} (when $f=1, A=1$). Not all values of the other parameters yield twisted tensor products; complete descriptions are provided at the beginning of the section in which a family is discussed.

\subsection{Noncommutative projective geometry} 
\label{nc proj geom}

We review some basic definitions and results from \cite{ATVI} and \cite{Mori}. Let $A = T(V)/\la R \ra$ be a quadratic algebra. An element $f\in T_d(V)$ determines a multilinear form $\tilde{f}:(V^*)^{\times d}\to \k$. Let $\G = \G(A) \subseteq \P(V^*) \times \P(V^*)$ denote the scheme of zeros of the multilinearizations of the elements of $R$. Let ${\rm pr}_i : \P(V^*) \times \P(V^*) \to \P(V^*)$ for $i = 1, 2$ denote the canonical projections. Let $E_i$ be the scheme-theoretic image of $\G$ under ${\rm pr}_i$. We denote ${\rm pr}_i|_{\G}$ by $\pi_i: \G \to E_i$.

\begin{defn}\label{def semi-standard}
We say that the algebra $A$ is \emph{semi-standard} if $E_1$ and $E_2$ are equal as subschemes of $\P(V^*)$.
\end{defn}

If the algebra $A$ is semi-standard, let $E = E_1 = E_2$. Then one may view $\G$ as the graph of a correspondence $E \to E$ via the closed immersion $(\pi_1, \pi_2): \G \to E \times E$. 

\begin{defn}\label{def nondegenerate}
Let $A$ be a semi-standard algebra. We say that $A$ is \emph{nondegenerate} if $\G$ is the graph of a scheme automorphism $\s: E \to E$. Otherwise we say that $A$ is \emph{degenerate}.
\end{defn}

If $A$ is semi-standard and nondegenerate, then we call the scheme $E$ the \emph{point scheme of $A$}. In \cite{ATVI} it is shown that the point scheme parametrizes the so-called point modules of $A$. We refer the reader to \cite{ATVI} for the definition of a point module of a graded algebra.

The pair $(E, \s)$ encodes important information about the algebra $A$. We use the notion of a \emph{geometric algebra} developed in \cite{Mori}. Let $i : E \to \P(V^*)$ be the natural embedding, and set $\L = i^* \O_{\P(V^*)}$. There is a natural map of $\k$-vector spaces: $$\m: H^0(E, \L) \tsr H^0(E, \L) \to H^0(E, \L \tsr_{\O_E} \sigma^* \L).$$

\begin{defn}\label{def geometric algebra}\cite[Definition 4.3]{Mori}
Let $A = T(V)/\la R \ra$ be a quadratic algebra. We say that $A$ is \emph{geometric} if there is a pair $(E, \s)$, where $E \subset \P(V^*)$ is a closed subscheme and $\s:E \to E$ is an automorphism such that:
\begin{itemize}
\item[(i)] $\G$ is the graph of $\s$,
\item[(ii)] $R = \ker \ \m$ under the canonical identification  $H^0(E, \L) = V$.
\end{itemize}
\end{defn}

We make frequent use of the following well-known result.

\begin{thm}\label{geometric algebra isomorphism} \cite[Remark 4.9]{Mori} Suppose that $A$ and $A'$ are geometric algebras associated to $(E, \s)$ and $(E', \s')$, respectively. Then $A \cong A'$ as graded $\k$-algebras if and only if $A_1 = A_1'$ and there is a scheme automorphism $\psi$ of $\P(A_1)$ which restricts to an isomorphism $\psi: E \to E'$ such that the following diagram commutes.
$$\xymatrix{&E \ar[r]^-{\psi} \ar[d]^(.445){\s} & E' \ar[d]^(.4){\s'} \\ &E \ar[r]^-{\psi} & E'
}
$$
\end{thm}

Some care must be taken in using the ``if direction'' of Theorem \ref{geometric algebra isomorphism} when the schemes $E, E'$ are not reduced. In the sequel, we do not use this part of the theorem in such cases.

In \cite{AS} the following notion of regularity for graded $\k$-algebras was introduced.

\begin{defn}\label{def AS regular}
\cite{AS}
A finitely-presented graded $\k$-algebra $A$ is \emph{Artin-Shelter regular (AS-regular) of dimension $d$} if:
\begin{itemize} 
\item[(i)] the global dimension of $A$ is $d$;
\item[(ii)] the Gelfand-Kirillov dimension of $A$ is finite;
\item[(iii)] $A$ is Gorenstein: $\Ext^n_A(\k, A) = 0$ if $n \ne d$, and $\Ext^d_A(\k, A) \cong \k$.
\end{itemize}
\end{defn}

Let us now restrict to the case we study in this paper and assume that $A = T(V)/\la R \ra$ is a quadratic algebra with $\dim_{\k} V =  \dim_{\k} R = 3$. In \cite{ATVI} it is shown that if $A$ is AS-regular, then $A$ is semi-standard and nondegenerate. Not all semi-standard, nondegenerate quadratic algebra are AS-regular, however. Following \cite{ATVI} we say that $A$ is \emph{exceptional} if the point scheme $E$ is the union of a line and a conic in $\P^2$ and the automorphism $\s: E \to E$ interchanges these components. In \cite[Proposition 4.11]{ATVI} it is proven that an exceptional algebra is not AS-regular. Finally, by \cite[Theorem 1, Theorem 6.8]{ATVI}, every quadratic AS-regular algebra of global dimension $3$ is geometric. 

Suppose that $R \subset V \tsr V$ has $\{r_1, r_2, r_3\}$ as a $\k$-basis. Using $x_i, y_i, z_i$, $i = 0, 1$, as coordinates on the $i+1$st factor of $\P(V^*) \times \P(V^*)$, we can factor the multilinearizations $\widetilde{r_1}, \widetilde{r_2}, \widetilde{r_3}$ in two ways:
\begin{equation}
\label{MandN}
M \begin{bmatrix} x_1 \\ y_1 \\ z_1 \end{bmatrix} = \begin{bmatrix} \widetilde{r_1} \\ \widetilde{r_2} \\ \widetilde{r_3} \end{bmatrix} , \qquad \begin{bmatrix} x_0 & y_0 & z_0 \end{bmatrix} N = \begin{bmatrix} \widetilde{r_1} & \widetilde{r_2} & \widetilde{r_3} \end{bmatrix},
\end{equation}
where $M$ is a $3 \times 3$ matrix with entries consisting of linear forms in $x_0, y_0, z_0$ and $N$ is a $3 \times 3$ matrix with entries consisting of linear forms in $x_1, y_1, z_1$. Then it is straightforward to prove that the set of closed points of the scheme $E_1$ is the zero locus of $\det M$, $\mathcal Z(\det M)$, in $\P^2$. Similarly, the set of closed points of $E_2$ is equal to $\mathcal Z(\det N)$ in $\P^2$. When $A$ is semi-standard, we often identify corresponding coordinate functions by suppressing the subscripts. It is easy to show that a semi-standard algebra is nondegenerate if and only if the rank of the matrix $M$ is $2$ at every closed point of $E$. Moreover, the automorphism $\s: E \to E$ can be computed on closed points by taking the cross product of two linearly independent rows of $M$; see \cite[Section 1]{ATVI}, for example.

For each AS-regular twisted tensor product in our classification, $A$, we determine the pair $(E, \s)$ for a particular representative of the graded isomorphism class of $A$. We adopt the terminology of \cite{MU} and \cite{IM} to describe the types of pairs that are possible. In \cite{MU} and \cite{IM} the ground field, $\k$, is assumed to be algebraically closed with ${\rm char } \, \k =0$. Though we assume only that $\k$ is algebraically closed and ${\rm char } \, \k \ne 2$, these  are sufficient for our purposes.

\begin{itemize}
\item Type ${\rm P}$: $E = \P^2$; $\s \in {\rm Aut}_{\k} \P^2 \cong {\rm PGL}_3(\k)$,
\item Type ${\rm S}_1$: $E$ is a triangle; $\s$ stabilizes each component,
\item Type ${\rm S}_2$: $E$ is a triangle; $\s$ interchanges two components,
\item Type ${\rm S}_3$: $E$ is a triangle; $\s$ cyclically permutes the components,
\item Type ${\rm S}_1'$: $E$ is a union of a line and a conic meeting at two points; $\s$ stabilizes each component and the intersection points,
\item Type ${\rm S}_2'$: $E$ is a union of a line and a conic meeting at two points; $\s$ stabilizes each component and interchanges the intersection points,
\item Type ${\rm T}_1$: $E$ is a union of three lines meeting at one point; $\s$ stabilizes each component,
\item Type ${\rm T}_2$: $E$ is a union of three lines meeting at one point; $\s$ interchanges two components,
\item Type ${\rm T}_3$: $E$ is a union of three lines meeting at one point; $\s$ cyclically permutes the components,
\item Type ${\rm T}'$: $E$ is a union of a line and a conic meeting at one point; $\s$ stabilizes each component,
\item Type ${\rm CC}$: $E$ is a cuspidal cubic curve,
\item Type ${\rm NC}$: $E$ is a nodal cubic curve,
\item Type ${\rm WL}$: $E$ is a union of a double line and a line,
\item Type ${\rm TL}$: $E$ is a triple line,
\item Type ${\rm EC}$: $E$ is an elliptic curve.
\end{itemize}

Note that some of the schemes described here are reducible and/or non-reduced. Furthermore, \cite{IM} divides these types into subtypes: Type P into Type ${\rm P}_i$ for $i = 1, 2, 3$; Type NC into Type ${\rm NC}_i$ for $i = 1, 2$; Type WL into Type ${\rm WL}_i$ for $i = 1, 2, 3$; and Type TL into Type ${\rm TL}_i$ for $i = 1, 2, 3, 4$.

In \cite[4.13]{ATVI}, the authors identified four subtypes of Type EC, denoted $A, B, E,$ and $H$, which were defined in terms of the automorphism $\s$. In \cite{Mat}, Matsuno completed the classification of algebras of Type EC up to graded algebra isomorphism in characteristic 0, showing that there are no additional subtypes.

\section{Elliptic-type twisted tensor products}
\label{elliptic type ttps}

Since we assume ${\rm char}\ \k\neq 2$,  \cite[Lemma 5.2]{CG3} implies that any elliptic-type twisted tensor product is isomorphic to an algebra in the family
$$T(g,h)=\dfrac{\k\la x,y,z\ra}{\la xy-yx, zy+yz-x^2-gy^2, z^2+hy^2\ra}\quad g, h\in \k.$$
In this section we show that all algebras in this family are semi-standard and nondegenerate. We describe the associated point schemes and describe the distinct graded algebra isomorphism classes in this family. 
Our results for elliptic-type twisted tensor products are summarized in Table \ref{EllipticTable} and the following theorem.

\begin{table}[h]
\caption{Mapping of elliptic-type TTP cases}
\label{EllipticTable}
\begin{tabular}{|c|c|c|c|c|}
\hline
TTP Case & Subcase & Algebra & Condition & Type\\
\hline
\multirow{4}{*}{$T(g,h)$} & \multirow{2}{*}{$h(g^2+4h)\neq 0$} & $T(g^2/h)$  &  $g^2/h \notin\{ 0, -4\}$   & \multirow{2}{*}{EC, subtype B} \\
\cline{3-4}
&&$T(0,1)$ &&\\
\cline{2-5}
& $h\neq 0=g^2+4h$ & $T(-4)$ &  & $NC_2$\\
\cline{2-5}
& $h= 0$ & $T(g,0)$ & $g\in\{0,1\}$ & exceptional\\
\hline
\end{tabular}
\end{table}

\begin{thm} 
\label{point scheme of T(g, h)}
Every elliptic-type twisted tensor product is isomorphic as a graded algebra to $T(g,h)$ for some $g,h\in\k.$ The algebras listed in Table \ref{EllipticTable} are pairwise non-isomorphic, and the table gives the types of AS-regular algebras.
\end{thm}

The proof of Theorem \ref{point scheme of T(g, h)} is at the end of this section.

In \cite{CG3}, we characterized when the algebras $T(g,h)$ are AS-regular. The first result in this section elaborates on that characterization, and describes the point scheme of $T(g,h)$ and the associated automorphism.

\begin{prop}
\label{T(g,h) is semistandard nondegenerate}
The algebra $T(g,h)$ is semi-standard and nondegenerate for any $g,h\in \k$. Moreover, $T(g,h)$ is AS-regular if and only if $h\neq 0$.
The point scheme of $T(g, h)$ is  $$E=\mathcal Z(hy^3+x^2z-yz^2+gy^2z)\subset \P^2.$$
The automorphism $\s : E \to E$ switches the points $[1:0:0]$ and $[0:0:1]$. If $h\neq 0$, $\s([x:y:z]) = [xz: yz: -hy^2]$ if $y \ne 0$. 

If $h=0$, $\s$ interchanges the line $z=0$ and the conic $x^2-yz+gy^2 = 0$.
\end{prop}

\begin{proof}
The matrices $M$ and $N$ defined in equation (\ref{MandN}) are $$M = \begin{bmatrix} -y_0 & x_0 & 0 \\ -x_0 & z_0 - gy_0 & y_0 \\ 0 & hy_0 & z_0 \end{bmatrix}, \ \ N = \begin{bmatrix} y_1 & -x_1 & 0\\ -x_1 & -gy_1+z_1 & hy_1 \\ 0 & y_1 & z_1 \end{bmatrix}.$$ If we identify coordinate functions by suppressing subscripts, one checks that $$\det M = hy^3+x^2z - yz^2+g y^2z = -\det N.$$ It follows that $T(g, h)$ is semi-standard. It is easy to prove that $\rank(M) \geq 2$ at all points of $\P^2$, so $T(g, h)$ is nondegenerate. By \cite[Theorem 6.2]{CG3}, the algebras $T(g,h)$ are AS-regular if and only if $h\neq 0$.

The descriptions of $\s([x:y:z])$ when $h\neq 0$ and $h=0$ are obtained by taking the cross product of appropriate rows in $M$.
\end{proof}

The algebras $T(g,0)$ are \emph{exceptional} in the sense of \cite[Section 4.9]{ATVI}. We prove that up to isomorphism, there are only two such algebras in the family $T(g,h)$.

\begin{prop}
\label{T(1,0) and T(0,0)}
If $g\neq 0$, $T(g,0)\cong T(1,0)$. Moreover, $T(1,0)\ncong T(0,0)$.
\end{prop}

\begin{proof}
The first statement follows by rescaling $z\mapsto gz$ and $x\mapsto \sqrt{g}x$. For the second statement, suppose there is an isomorphism $\varphi:T(1,0)\to T(0,0)$. It is easy to check that, up to rescaling, $z$ is the only degree-1 element of $T(0,0)$ that squares to 0. Thus we may assume $\varphi(z)=z$.  

It is also straightforward to check that if $u,v\in T(0,0)_1$ are linearly independent and commute, then $u,v\in{\rm span}_{\k}\{x,y\}$. Thus $\varphi(x)=\a_1x+\b_1y$ and $\varphi(y)=\a_2x+\b_2y$ for $\a_1,\a_2,\b_1,\b_2\in\k$. Now, a direct calculation shows that $\varphi(zy+yz-x^2-y^2)=0$ implies $\varphi(x)$ and $\varphi(y)$ are linearly dependent, a contradiction.
\end{proof}

For the remainder of the section, we assume $h\neq 0$. As shown above, all algebras $T(g,h)$ are AS-regular in this case. 

When $h\neq 0$, the closed points of $E$ describe the projectivization of a plane cubic curve. After the transformation $x\mapsto y$, $y\mapsto -x$, $z\mapsto hz$, we see that $E$ can be described in Weierstrass form by $$y^2z=x^3-gx^2z-hxz^2.$$
Following the formulae in \cite[III.1]{Sil}, the discriminant of the curve obtained by setting $z=1$ is $\Delta(g,h) = 16h^2(g^2+4h)$. Thus we see that $E$ is the projectivization of an elliptic curve if and only if $g^2+4h\neq 0$. Furthermore, the $j$-invariant of this elliptic curve is $j(g,h) = \dfrac{16^2(g^2+3h)^3}{h^2(g^2+4h)}$.
Setting $\ell=g^2/h$ yields the expression $j(\ell)=\dfrac{16^2(\ell+3)^3}{\ell+4}$, which is sometimes more convenient.
Recall that the $j$-invariant classifies elliptic curves up to isomorphism, see \cite[Proposition III.1.4 (b)]{Sil}.

The next result characterizes $T(g,h)$ when $E$ is a singular curve; the case where $E$ is an elliptic curve is described below.

\begin{prop}
\label{nodal cubic}
If $h\neq 0$ and $g^2+4h=0$, then $E$ is a nodal cubic curve, and the automorphism $\s$ fixes the node. (type $NC_2$)\\
In this case, $T(g,h)$ is isomorphic as a $\k$-algebra to
$$\dfrac{\k\la x,y,z\ra}{\la xz-2yx+zy, zx-2xy+yz, y^2+x^2\ra}.$$
\end{prop}

\begin{proof}
Using the description of $E$ as $\mathcal Z(hy^3+x^2z-yz^2+gy^2z)$ from Theorem \ref{T(g,h) is semistandard nondegenerate}
 we see that 
if $h\neq 0$ and $g^2+4h=0$, then $E$ is the projectivization of a nodal cubic curve: $$x^2z=\dfrac{y(2z-gy)^2}{4}$$
Evidently the node is $[0:2:g]$, and by Theorem \ref{T(g,h) is semistandard nondegenerate} we see that $\s([0:2:g])=[0:2g:-4h]=[0:2:g]$ since $g\neq 0$.

Let $A$ denote the algebra specified in the statement of the Proposition. The map determined by: $x\mapsto \a_1(x-y)$, $y\mapsto \a_2(x+y+z)$, $z\mapsto -3x-3y+z$ where $\a_2=-g/2h$ and $\a_1^2=-8g/h$, is an isomorphism $T(g, h) \to A$.
\end{proof}

It remains to consider the case where the discriminant $\Delta(g,h)\neq 0$. Our next result shows that in this case, $T(g, h)$ is a Type EC, subtype B algebra.

\begin{thm}
\label{sigma is multiplication by -1} 
Suppose that $\Delta(g,h) \ne 0$, so that $E$ is a nonsingular cubic curve. There exists an element $O \in E$ such that for the corresponding group law on $E$, the automorphism $\s: E \to E$ is given by $\s(P)+P=O$ for all $P\in E$.
\end{thm}

\begin{proof}
It is straightforward to check that a point $[x_0:1:z_0]\in E$ is a fixed point of $\s$ if and only if $x_0^2=2z_0-g$ and $z_0^2=-h$. Let $O\in E$ be any fixed point of $\s$. Consider the corresponding group law on $E$, see \cite[III.2]{Sil} for example.

The line $y=0$ intersects $E$ in $e_1=[1:0:0]$ and $e_3=[0:0:1]$, and $\s$ interchanges these points. We claim that $e_1+e_3=O$ in the group law.
To see this, first observe that $y=0$ is tangent to $E$ at $e_3$. The line $L$ through $e_3$ and $O=[x_0:1:z_0]$ is $x=x_0y$. The point $e_1$ is not on $L$, and $L$ is not tangent to $E$ at $e_3$, so we can write $e_1+e_3=[x_0:1:z_1]$. The condition that $e_1+e_3$ lies on $\mathcal Z(hy^3+x^2z-yz^2+gy^2z)$, and the fact that $x_0^2=2z_0-g$ and $z_0^2=-h$ implies
$$ z_1^2-2z_0z_1+z_0^2 = (z_1-z_0)^2 = 0.$$
Thus $z_1=z_0$ and $e_1+e_3=0$ as desired. 

The argument above also proves that $x=x_0y$ is the tangent line to $E$ at $O$. Since we assumed only that $O$ was a fixed point of $\s$ on $E$, if $[x_1:1:z_1] \in E$ is any other fixed point of $\s$, the tangent line at that point is $x=x_1y$. 

Let $P = [x_1: 1: z_1] \in E-\{e_1,e_3\}$ be arbitrary. Then $\s(P) = [x_1 z_1: z_1: -h]$. If $\s(P)\neq P$, then $z_1^2 \ne -h$.  The line through $P$ and $\s(P)$ is  $x - x_1 y=0$, and $e_3$ is the third point of intersection between this line and $E$. The line through $e_3$ and $O$ is tangent to $E$ at $O$, so $P + \s(P) = O$.

If $\s(P)=P$, then as observed above, the tangent line to $E$ at $P$ is the line $x -x_1 y=0$. The calculation then proceeds as before, yielding $P+\s(P) = O$. 
\end{proof}

Next we classify the algebras $T(g,h)$ with $h\neq 0$ up to graded algebra isomorphism.

\begin{prop}
\label{isomorphism type of T(g,h)}
If $h\neq 0$ and $h' \ne 0$, then $T(g,h)\cong T(g',h')$ if and only if $g^2h'=(g')^2h$. 
\end{prop}

\begin{proof}
Assume that $g^2h'=(g')^2h$. If $g \ne 0$, then the map determined by $x \mapsto \a x$, $y \mapsto \b y$, $z \mapsto z$, where $\a^2 = \b$ and $\b = g'/g$ is an isomorphism $T(g, h) \to T(g', h')$. If $g = 0$, then the map determined by $x \mapsto \a x$, $y \mapsto \b y$, $z \mapsto z$, where $\a^2 = \b$ and $\b^2 = h'/h$ is an isomorphism $T(g, h) \to T(g', h')$.

Suppose, conversely, that $\varphi:T(g',h')\to T(g,h)$ is an isomorphism.
A calculation, using the bases $\{x^2, xy, xz, y^2, yz, zx\}$ in degree $2$, shows that commuting, linearly independent degree-1 elements of $T(g,h)$ lie in ${\rm span}_{\k}\{x,y\}$. Thus there is no loss of generality in assuming that
$$
\varphi(x)=\a_1x+\b_1y,\quad
\varphi(y)=\a_2x+\b_2y,\quad
\varphi(z)=\a_3x+\b_3y+z.
$$
Then 
\begin{align*}
0=\varphi(z^2+h'y^2)&=(\a_3^2+\b_3+h'\a_2^2)x^2+2(\a_3\b_3+h'\a_2\b_2)xy\\
&\quad +(\b_3^2+g\b_3-h+h'\b_2^2)y^2+\a_3(xz+zx),
\end{align*}
so we have
$$\a_3=0\quad \b_3=-h'\a_2^2\quad \a_2\b_2=0\quad \b_3^2+g\b_3+h'\b_2^2-h=0.$$
Since $\a_3=0$, the coefficient of $xz$ in $\varphi(zy+yz-x^2-g'y^2)$ is $\a_2$, so $\a_2=0$, $\b_3=0$, and $h'\b_2^2=h$. Bearing these in mind, we have
$$\varphi(zy+yz-x^2-g'y^2)=(\b_2-\a_1^2)x^2-2\a_1\b_1xy+(g\b_2-g'\b_2^2-\b_1^2)y^2.$$
We cannot have $\a_1=0$, else $\b_2=0$ and $\varphi(y)=0$, so $\b_1=0$ and $g\b_2=g'\b_2^2$. Since $\b_2\neq 0$, we have
$h'\b_2^2=h$ and $g'\b_2=g$, and the result follows.
\end{proof}

\begin{rmk}
\label{T(g,h) remark}
Proposition \ref{isomorphism type of T(g,h)} motivates the following further simplification. For $\ell\in \k^*$, we define
$$T(\ell)=\k\la x,y,z\ra/\la xy-yx, zy+yz-x^2-\ell y^2, z^2+\ell y^2\ra.$$
When $gh\neq 0$, it is easily shown that $T(g,h)\cong T(g^2/h)$ by rescaling $x\mapsto \sqrt{g/h}x$ and $y\mapsto gy/h$. Proposition \ref{isomorphism type of T(g,h)} then implies $T(\ell)\cong T(\ell')$ if and only if $\ell=\ell'$. 

Recall that the $j$-invariant of the point scheme of $T(\ell)$ is given by $j(\ell)=\dfrac{16^2(\ell+3)^3}{\ell+4}$. Setting $l = 0$ we see that $j(0) = 12^3$. The point scheme of $T(0,1)$ has $j=12^3$; we ``compactify" the family of algebras $T(\ell)$, $\ell \in \k^*$ by defining $T(0) = T(0,1)$. (The algebra obtained by setting $\ell=0$ in the presentation above is $T(0,0)$, which is not AS-regular.)

The $j$-invariant does not completely characterize the isomorphism type of the algebra $T(\ell)$. Setting $t=\ell+3$ and considering the discriminant of the polynomial $t^3-(j/16^2)(t+1),$
we see that for a fixed elliptic curve $E$, or equivalently for a fixed value of $j$, there are three isomorphism classes of algebras $T(g,h)$ with point scheme equal to $E$ whenever $j\notin \{0, 12^3\}$, two when $j=12^3$, and one when $j=0$. 

Comparing the above to \cite[Proposition 3.10]{Mat} shows that in characteristic 0, every quadratic AS-regular algebra of type EC, subtype B, is isomorphic to a graded twisted tensor product of $\k[x,y]$ and $\k[z]$. By \cite[Theorem 3.1]{IM}, there is only one isomorphism class of type $NC_2$ algebras, and Proposition \ref{nodal cubic} shows these algebras are graded twisted tensor products of $\k[x,y]$ and $\k[z]$.
\end{rmk}

We conclude with a proof of Theorem \ref{point scheme of T(g, h)}.

\begin{proof}[Proof of Theorem \ref{point scheme of T(g, h)}]
Every elliptic-type twisted tensor product is isomorphic to $T(g,h)$ by \cite[Lemma 5.2]{CG3}. Proposition \ref{T(g,h) is semistandard nondegenerate} shows that $T(g,h)$ is AS-regular if and only if $h\neq 0$. 
The algebras $T(g,0)$ are classified  in Proposition \ref{T(1,0) and T(0,0)}.

When $h\neq 0$, the graded algebra isomorphism classes of $T(g,h)$ are described in Proposition \ref{isomorphism type of T(g,h)}. This classification is described in terms of the family $T(\ell)$, $\ell\neq 0$, in Remark \ref{T(g,h) remark}. The types of AS-regular algebras in the family $T(g,h)$ are characterized in Proposition \ref{nodal cubic} and Theorem \ref{sigma is multiplication by -1}.
\end{proof}

\section{Reducible-type twisted tensor products}
\label{reducible type ttps}

In this section, $T$ denotes a quadratic twisted tensor product of reducible type, in \emph{normal form} \cite[Section 4]{CG3}. The results below classify $T$ up to graded algebra isomorphism. When $T$ is AS-regular, we identify its type. We assume  ${\rm char}\ \k\neq 2$. 

By \cite[Theorem 4.9]{CG3} $T\cong \k\la x,y,z\ra/I$ where $I$ is the graded ideal of $\k\la x,y,z\ra$ generated by the elements
$$xy-yx, zx-ax^2-bxy-cy^2-dxz-eyz-z^2, zy-Bxy-Cy^2-Dyz,$$
such that $(a, d)$ is not a zero of certain polynomials $f_n(t,u)$ for all $n \geq 1$, and the parameters satisfy one of the following cases:
\begin{itemize}
\item[R(i):] $a = B(1-d-B)$, $b = c = e=C = D = 0$;
\item[R(ii):] $e=B = C = 0$, $D = 1$;
\item[R(iii):] $d = -1$, $B = 2$, $e=C = 0$, $D = -1$;
\item[R(iv):] $a = B(1-d-B)$, $b = 1-d-2B$, $c=-1$, $C = 1$, $e=D = 0$;
\item[R(v):] $e=0$, $d= D= -1$, $B = 2$, $C = 1$.
\item[R(vi):] $e=1$, $d=D=0$, $a=B(1-B)$, $b=C-B-2BC$, $c=-C(1+C)$;
\item[R(vii):] $B=C=0$, $e=d=D=1$.
\end{itemize}

We use the polynomials $f_n(t,u)$ in the proof of Proposition \ref{non-regular classification} below. We note here that $f_1(t,u)=1-t$, so $a \ne 1$.

By \cite[Theorem 6.2]{CG3}, no algebras in cases R(i), R(iv), R(vi) are AS-regular. After two preliminary results, the analyses in this section are divided into two subsections, according to whether or not $T$ is AS-regular.
The results are summarized in Table \ref{ReducibleTable} and the following theorem. 

\begin{table}[h]
\caption{Mapping of reducible-type TTP cases not isomorphic to an Ore-type TTP}
\label{ReducibleTable}
\setlength\tabcolsep{2pt}
\begin{tabular}{|c|c|c|c|}
\hline
TTP Case & Subcase & Algebra   & Type\\
\hline
 \multirow{2}{*}{R(iii)} & $b^2-4ac+4c=0$ & $U$ & $S_2$\\
 \cline{2-4}
 &$b^2-4ac+4c\neq 0$&  $U'$ &  $S_2'$\\
 \hline
 \multirow{2}{*}{R(i)} & $B+d=0$  &  $R_0$ & \multirow{4}{*}{not AS-regular} \\
  \cline{2-3}
 & \multirow{2}{*}{$B+d\neq 0$} & \multirow{2}{*}{$R(q)$} &  \\
 \cline{1-1}
 \multirow{2}{*}{R(iv)} & & &\\
 \cline{2-3}
  & $B+d= 0$  &  $R_1$ &  \\
\hline
\end{tabular}
\end{table}

\begin{thm}
\label{ReducibleSummary}
Every reducible-type twisted tensor product is isomorphic as a graded algebra to one of the algebras listed in Table \ref{ReducibleTable}, or to a twisted tensor product of Ore type. The algebras listed in Table \ref{ReducibleTable} are pairwise non-isomorphic, and the table gives the types of the AS-regular algebras.
\end{thm}

The proof of Theorem \ref{ReducibleSummary} can be found at the end of the section.

\begin{rmk}
\label{complete reducible cases}
By \cite[Theorem 3.1]{IM}, in characteristic 0, the isomorphism classes of Type $S_2$ AS-regular algebras are parametrized by points of $\P^1$, but there is only one isomorphism class of Type $S'_2$. We show in Section 5 that no other graded twisted tensor products of $\k[x,y]$ and $\k[z]$ are of Type $S_2$.
\end{rmk}

The first lemma, and the theorem that follows, show that only case R(iii) contains AS-regular algebras not isomorphic to a graded Ore extension of $\k[x,y]$. (Ore-type twisted tensor products are discussed in Section \ref{ore type ttps}.)
\begin{lemma}
\label{reduction to case (iii)} The following hold.

\begin{itemize}
\item[(1)] Any twisted tensor product from case R(v) is isomorphic as a graded algebra to an algebra from case R(iii).
\item[(2)] Any twisted tensor product from case R(vii) is isomorphic to an algebra from case R(ii).
\item[(3)] Any twisted tensor product from case R(ii) is isomorphic to an Ore-type twisted tensor product.
\end{itemize}
\end{lemma}

\begin{proof} 
For (1), the map sending $x\mapsto x-y/2$ and fixing $y$ and $z$ identifies a case R(v) algebra with one from case R(iii).

For (2), the map sending $z\mapsto z-y/2$ and fixing $x$ and $y$ identifies a case R(vii) algebra with one from case R(ii).

For (3), If $T$ is a reducible-type twisted tensor product from case R(ii) with $a=0$ and $d=1$, then the map interchanging $x$ and $z$, fixing $y$, identifies $T$ with the Ore-type algebra with relations
$$xy-yx, zy-yz, zx+x^2+cy^2-xz-byz.$$
If $T$ is any other reducible-type twisted tensor product from case R(ii), the map sending $x\to x+\a z$ and fixing $y$ and $z$ where $\a$ satisfies $a\a^2+(d-1)\a+1=0$ identifies a case R(ii) algebra with the algebra with relations
$$xy-yx, zy-yz, (1-a\a)zx-ax^2-bxy-cy^2-(d+a\a)xz-b\a yz.$$
Note that this algebra is an Ore-type twisted tensor product provided that $1-a\a \ne 0$. If $a = 0$, then this is clear.  Suppose that $a \ne 0$. Let $\d = (d-1)^2-4a$ be the discriminant of $at^2+(d-1)t+1$. If $\d \ne 0$, then $a\a^2+(d-1)\a+1=0$ has two distinct solutions and we may choose $\a \ne a^{-1}$. Suppose that $\d = 0$ and $\a = a^{-1}$ is the unique solution of $a\a^2+(d-1)\a+1=0$. One checks that $a = 1$, but this contradicts the fact that $T$ is a twisted tensor product. 
\end{proof}

\begin{thm}
\label{semistandard reducible type}
Let $T$ be a reducible-type twisted tensor product in normal form. Consider the following statements.
\begin{enumerate}
\item $a+d\neq 0\neq D$,
\item $T$ is AS-regular of dimension 3,
\item $T$ is semi-standard.
\end{enumerate}
Then $(1)\Rightarrow (2)\Rightarrow (3) \Rightarrow D\neq 0$. If $T$ is semi-standard and $a+d=0$, then $T$ is isomorphic as a graded algebra to an Ore-type twisted tensor product.
 \end{thm}

\begin{proof}
The implication $(1)\Rightarrow (2)$ follows from \cite[Theorem 6.2(2)]{CG3}, and  $(2)\Rightarrow (3)$ is immediate. Assume $T$ is semi-standard. The matrices $M$ and $N$ of equation (\ref{MandN}) for a reducible-type algebra are
$$M=\begin{bmatrix}
-y_0 & x_0 & 0\\
z_0-ax_0 & -bx_0-cy_0 & -dx_0-ey_0-z_0\\
0 & z_0-Bx_0-Cy_0 & -Dy_0
\end{bmatrix}$$
and
$$N=\begin{bmatrix}
y_1 & -ax_1-by_1-dz_1 & -By_1\\
-x_1 & -cy_1-ez_1 & -Cy_1-Dz_1\\
0 & x_1-z_1 & y_1
\end{bmatrix}.
$$ 
Suppressing subscripts we have
\begin{align*}
\det M &= -y( (aD-dB)x^2+(bD-eB-dC)xy+(cD-eC)y^2\\
&\quad +(d-B-D)xz+(e-C)yz+z^2),\\
\det N &= y( (B-a)x^2 +(C-b)xy -cy^2 +(D-d-B)xz -(e+C)yz -Dz^2).
\end{align*}
Since $T$ is semi-standard, we have $\det M=k\det N$ for some $k\in \k$. Examining the coefficients of $yz^2$, we see that $D\neq 0$, so $T$ belongs to case R(ii), R(iii), R(v), or R(vii). 
If $T$ is an algebra from case R(iii) or R(v), then $d=-1$. Since $f_1(a,d)=1-a\neq 0$, it follows that $a+d\neq 0$. 
If $T$ is  from case R(ii) or R(vii), then $T$ is isomorphic to an Ore-type twisted tensor product by Lemma \ref{reduction to case (iii)}.
\end{proof}

\subsection{Twisted tensor products from case R(iii)}

In case R(iii), $a+d\neq 0$, and Theorem \ref{semistandard reducible type} implies that 
all algebras in this family are AS-regular. If $T$ is an algebra from case R(iii), the matrix $M$ is
$$\begin{bmatrix}
-y_0 & x_0 & 0\\
z_0-ax_0& -bx_0-cy_0 & x_0-z_0\\
0 & z_0-2x_0 & y_0
\end{bmatrix},$$
so the point scheme of $T$ is
$$E=\mathcal Z(y(z^2-2xz+(2-a)x^2-bxy-cy^2)).$$
The automorphism, which has order 2, is
$$\s([x:y:z])=\begin{cases} [x:y:2x-z] & y\neq 0\\ [z-x:0:z-ax] & y=0\\ \end{cases}.$$
As ${\rm char}\ \k\neq 2$, the quadratic form $z^2-2xz+(2-a)x^2-bxy-cy^2$ is reducible if and only if its discriminant, $b^2-4ac+4c=0$.

Despite the number of parameters, there are only two graded algebra isomorphism classes of case R(iii) algebras, according to whether or not $b^2-4ac+4c=0$, as the next result shows. This result also clarifies the action of $\s$ on $E$. We define the two algebras as follows. Let
$$U=\k\la x,y,z\ra/\la xy-zx, yx-xz, y^2+z^2\ra,$$ 
and
$$U'=\k\la x,y,z\ra/\la xy-zx, yx-xz, x^2+y^2+z^2\ra.$$

\begin{prop}
\label{S_2',S_2}
Let $T$ be a twisted tensor product from case R(iii). Let $s = b^2-4ac+4c$.
\begin{itemize}
\item[(1)] If $s \ne 0$, then $T$ is isomorphic to the algebra $U'$.
\item[(2)] If $s = 0$, then $T$ is isomorphic to the algebra $U$.
\end{itemize}
\end{prop}

\begin{proof}
For (1), let  $\a_2=2(1-a)\sqrt{-2/s}$ and let $\a_1=\dfrac{b\a_2}{2(1-a)}$. Since $a\neq 1$, we have $\a_2\neq 0$. Define a graded algebra homomorphism $\varphi:\k\la x,y,z\ra\to \k\la x,y,z\ra$ by
\begin{align*}
\varphi(x) &= \a_1x+y+z\\
\varphi(y) &= \a_2x\\
\varphi(z) &= \a_1x+(1+\sqrt{a-1})y+(1-\sqrt{a-1})z.
\end{align*}
Since $\a_2\neq 0$, and $a-1\neq 0$, the map $\varphi$ is clearly bijective in degree 1, and hence is an isomorphism. It suffices to prove that that $\varphi$ is a $\k$-linear isomorphism from the quadratic relations space of $T$ to that of $U'$. For ease of notation we define $\b_3=1+\sqrt{a-1}$ and $\d_3=1-\sqrt{a-1}$. By direct calculation we see that
$$\varphi(xy-yx)=\a_2(yx+zx-xy-xz)$$
and 
$$\varphi(zy-2xy+yz)=\a_2\left[(\b_3-2)yx-\b_3zx+\b_3xy-(\b_3-2)xz\right].$$
Since $\b_3\neq 1$, these two images span the same subspace as $yx-xz$ and $zx-xy$. We consider $\varphi(zx-ax^2-bxy-cy^2+xz-z^2)$ modulo these two relations, and show that it is a nonzero scalar multiple of $x^2+y^2+z^2$. To that end,

\begin{align*}
\varphi(zx-z^2) &= (1-\b_3)(\a_1x+\b_3y+\d_3z)(y-z)\\
&=(1-\b_3)(\a_1xy-\a_1xz+\b_3y^2-\b_3yz+\d_3zy-\d_3z^2),\\
\varphi(-ax^2) &= -a(\a_1^2x^2+2\a_1xy + 2\a_1xz + y^2+yz+zy+z^2),\\
\varphi(-bxy) &= -b\a_2( \a_1x^2+xz+xy), \\
\varphi(-cy^2) &= -c\a_2^2x^2, \\
\varphi(xz) &= ( \a_1x+y+z)(\a_1x+\b_3y+\d_3z),\\
&=\a_1^2x^2+\a_1(1+\b_3)xy+\a_1(1+\d_3)xz+\b_3y^2+\d_3yz\\
&\quad +\b_3zy+\d_3z^2.
\end{align*}
Examining the coefficients of the sum of the terms above, and noting that $1-\b_3=\d_3-1$, we see that the coefficients of $xy$ and $xz$ both equal $(2-2a)\a_1-b\a_2$, which vanishes by definition of $\a_1$. The coefficients of $yz$ and $zy$ both equal $\b_3^2-2\b_3+2-a$, which vanishes by definition of $\b_3$. The coefficients of $y^2$ and $z^2$ both equal $2-2a$, and the coefficient of $x^2$ is $(1-a)\a_1^2-b\a_2\a_1-c\a_2^2.$ Using the definitions of $\a_1$ and $\a_2$, it is straightforward to show this expression also equals $2-2a$.

The proof of (2) is similar. Let $\d_1=\dfrac{b}{2(1-a)}$. Define a graded algebra homomorphism $\varphi:\k\la x,y,z\ra\to \k\la x,y,z\ra$ by
\begin{align*}
\varphi(x) &= \d_1x+y+z\\
\varphi(y) &= x\\
\varphi(z) &= \d_1z+(1+\sqrt{a-1})x+(1-\sqrt{a-1})y.
\end{align*}
One then checks that this map induces the desired isomorphism $T\to U$. The details are left to the reader. 
\end{proof}

In light of Lemma \ref{reduction to case (iii)} and Proposition \ref{S_2',S_2}, there are only two isomorphism classes of AS-regular twisted tensor products of reducible type that do not contain Ore-type twisted tensor products. The next result shows these isomorphism classes are distinct, and characterizes their types.

\begin{prop}
\label{reduced-type case (iii)}
The algebras $U$ and $U'$ are AS-regular. 
\begin{itemize}
\item[(1)] The point scheme of $U'$ is a union of a line and an irreducible conic, meeting at two points. The automorphism $\s$ stabilizes the components and interchanges the intersection points. (type $S_2'$)
\item[(2)] The point scheme of $U$ is a union of three lines meeting in three points, and $\s$ interchanges two of the lines, fixing $y=0$. (type $S_2$)
\end{itemize}
In either case, the automorphism $\s$ has order 2.
\end{prop}

\begin{proof}
By Proposition \ref{S_2',S_2}, both $U$ and $U'$ are isomorphic to case R(iii) algebras, which are all AS-regular by Theorem \ref{semistandard reducible type}.

First consider the algebra $U'$. The associated matrix $M$ of equation (\ref{MandN}) is
$$
\begin{bmatrix}
-z_0 & x_0 & 0\\
y_0 & 0 & -x_0\\
x_0& y_0& z_0
\end{bmatrix},
$$
so the point scheme is $E=\mathcal Z(x(2yz+x^2))$, a union of a line and an irreducible conic. The line $x=0$ intersects the conic $2yz+x^2 = 0$ in $[0:1:0]$ and $[0:0:1]$. The automorphism $\s$ is given by $$\s([x:y:z]) = \begin{cases} [x: z: y] & x\neq 0\\ [0: -z: y] & x=0.\\ \end{cases}$$
So $\s$ stabilizes the components of $E$, interchanges $[0:1:0]$ and $[0:0:1]$, and has order 2.

Next consider the algebra $U$. The matrix $M$ of equation (\ref{MandN}) is
$$
\begin{bmatrix}
-z_0 & x_0 & 0\\
y_0 & 0 & -x_0\\
0 & y_0& z_0
\end{bmatrix},
$$
thus the point scheme is $E=\mathcal Z(xyz)$, a union of three lines. The automorphism $\s$ is given by
\begin{align*}
\s([0:y:z]) &= [0: -z: y]\\
\s([x:0:z]) &= [x: z: 0]\\
\s([x:y:0]) &= [x:0:y],
\end{align*}
so $\s$ stabilizes $x=0$ while interchanging $y=0$ and $z=0$. The automorphism $\s$ has order 2.
\end{proof}

\subsection{Twisted tensor products from cases R(i), R(iv), R(vi)}

Recall that, as noted above, none of the reducible-type algebras in cases R(i), R(iv) and R(vi) are AS-regular. It appears that the algebras in each respective case may depend, even up to isomorphism, on two parameters. However, we show that, up to isomorphism, the algebras in cases R(i), R(iv) and R(vi) constitute a 1-parameter family of algebras and two isolated algebras.

\begin{lemma}
\label{case (vi) is (i) or (iv)}
Any twisted tensor product from case R(vi) is isomorphic to an algebra from case R(i) or case R(iv).
\end{lemma}

\begin{proof}
The defining relations of a twisted tensor product $T$ from case R(vi) are
$$xy-yx, zy-Bxy-Cy^2, zx-ax^2-bxy-cy^2-yz-z^2,$$ where $a = B(1-B)$, $b = C-B-2BC$, $c = -C(1+C)$. 
The map fixing $x$ and $y$, sending $z\mapsto z+Cy$ is an isomorphism from $T$ to the algebra $T'$ with relations
$$xy-yx, zy-Bxy, zx-ax^2+B(1+C)xy-(1+C)yz-z^2.$$
If $C=-1$, this algebra belongs to case R(i) ($d = 0$). Otherwise, the map sending $y\mapsto y/(1+C)$, $z\mapsto z-y$ and fixing $x$ is an isomorphism from $T'$ to an algebra from case R(iv) ($d = 0$).
\end{proof}

For $d, B, q \in \k$ we define 
\begin{align*}
&T_{(i)}(d,B)=\k\la x,y,z\ra/\la xy-yx, zy, (1-B)zx-(d+B)xz-z^2\ra, \\
&T_{(iv)}(d,B)=\k\la x,y,z\ra/\la xy-yx, zy, (1-B)zx-(d+B)xz-yz-z^2\ra, \\
&R(q)=\k\la x,y,z\ra/\la xy-yx, zy, qzx-xz-z^2\ra, \\
&R_0 = \k\la x,y,z\ra/\la xy-yx, zy, zx-z^2\ra, \\
&R_1 = \k\la x,y,z\ra/\la xy-yx, zy, zx-yz-z^2\ra. 
\end{align*}

We omit the straightforward verification of the following lemma.

\begin{lemma}
\label{R(q) are distinct}
As graded algebras, $R(q)\cong R(q')$ if and only if $q=q'$. Furthermore, for all $q \in \k$, the algebras $R(q)$, $R_0$ and $R_1$ are pairwise non-isomorphic. 
\end{lemma}

We note that for any $\l \in \k$, $\l \ne 0$, $R_1 \cong \k\la x,y,z\ra/\la xy-yx, zy,zx-\l yz-z^2\ra,$ so one could also view $R_0$ and $R_1$ as members of a 1-parameter flat family.

\begin{prop} 
\label{non-regular classification}
An algebra $T$ is a twisted tensor product from case R(i) or R(iv) if and only if it is isomorphic as a graded algebra to one of $R(q)$, $q \in \k$, $R_0$, or $R_1$. 
\end{prop}

\begin{proof}
First, let $T$ be a twisted tensor product from case R(i) or case R(iv).  If $T$ is from case R(i), then the map that fixes $x$ and $y$ and sends $z \mapsto Bx+z$ is an isomorphism $T \to T_{(i)}(d,B)$. If $T$ is from case R(iv), then the map that fixes $x$ and $y$ and sends $z \mapsto Bx+y+z$ is an isomorphism $T \to T_{(iv)}(d,B)$.

If $B+d\neq 0$, the map given by $x\mapsto x+y$, $y\mapsto -(d+B)y$, $z\mapsto z$ is an isomorphism $T_{(iv)}(d,B)\rightarrow T_{(i)}(d,B)$. Rescaling $x$ by $1/(d+B)$ and fixing $y$ and $z$ gives an isomorphism between $T_{(i)}$ and $R(q)$, where $q = (1-B)/(B+d)$.

If $B+d=0$, the requirement that $a=B(1-d-B)$ implies $a=B$, and since $a \ne 1$, it follows that $B\neq 1$. Then one sees that $T_{(i)}(d,B) \cong R_0$ and $T_{(iv)}(d,B) \cong R_1$  by isomorphisms fixing $y$ and $z$ and sending $x \mapsto (1/(1-B))x$.

For the converse, we first note that by \cite[Lemma 3.3]{CG3}, $f_n(a,-a)=(1-a)^n$ for all $n\ge 0$. Thus choosing $a=-d\neq 1$ guarantees $(a,d)$ is not a zero of any $f_n(t,u)$. Now it is clear that $R_0$ and $R_1$ are isomorphic to case R(i) and R(iv) twisted tensor products where, for example, $a=B=-d=0$. Similarly, $R(0)$ is isomorphic to a twisted tensor product from case R(i)  where $B=1$ and $a=d=0$. 

If $q\in \k^*$, put $d=1/q$ and $a=B=0$. One can show that in this case $f_n(a,d)=1$ for all $n\ge 0$, so there exists a case R(i) twisted tensor product isomorphic to $R(q)$. 
\end{proof}

We conclude this section with the proof of Theorem \ref{ReducibleSummary}.

\begin{proof}[Proof of Theorem \ref{ReducibleSummary}]
By Lemma \ref{reduction to case (iii)}(2,3), the algebras in cases R(ii) and R(vii) are isomorphic to Ore-type twisted tensor products. Lemma \ref{reduction to case (iii)}(1) and Lemma \ref{case (vi) is (i) or (iv)} imply that among the remaining cases, it suffices to classify R(i), R(iii), and R(iv). 

Proposition \ref{S_2',S_2} shows that every algebra in case R(iii) is isomorphic to $U$ or $U'$. Proposition \ref{reduced-type case (iii)} shows these algebras are non-isomorphic and AS-regular, and identifies their types as $S_2$ and $S_2'$, respectively.

By \cite[Theorem 6.2]{CG3}, no algebras in case R(i) or R(iv) are AS-regular. These algebras are classified up to graded algebra isomorphism as $R(q)$, $R_0$, or $R_1$ in Proposition \ref{non-regular classification}, and are pairwise non-isomorphic by Lemma \ref{R(q) are distinct}.
\end{proof}

\section{Ore-type twisted tensor products}
\label{ore type ttps}

In this section, let $T$ denote a quadratic twisted tensor product of Ore type, in normal form.  The results below classify $T$ up to isomorphism of graded algebras. If $T$ is AS-regular, we identify its type, as in \cite{IM}. We continue to assume ${\rm char}\ \k\neq 2$.

We begin by recalling the description of Ore-type twisted tensor products from \cite{CG3}. By \cite[Theorem 4.3]{CG3}, $T\cong \k\la x,y,z\ra/I$ where $I$ is the graded ideal of $\k\la x,y,z\ra$ generated by the elements
$$xy-yx, zx-ax^2-bxy-cy^2-dxz-eyz, zy-Ax^2-Bxy-Cy^2-Dyz,$$
such that $A, e\in \{0,1\}$ and $a, b, c, d, B, C, D\in \k$ belong to one of the following cases. 

\begin{enumerate}
\item[O(1):] $e=0$, $A\in \{0,1\}$, and if $A=0$, then $C\in \{0,1\}$. Furthermore, one of the following holds:
\begin{enumerate}
\item[(i)] $d=D=1$,
\item[(ii)]  $d\neq 1=D$, $A=B=C=0$,
\item[(iii)] $d=1\neq D$, $a=b=c=0$,
\item[(iv)] $d\neq 1\neq D$, $A=c=0$, $a=B(d-1)/(D-1)$, $b=C(d-1)/(D-1)$.
\end{enumerate}
\vspace{0.25cm}
\item[O(2):] $e=1$, $d=D$ and one of the following holds:
\begin{enumerate}
\item[(i)] $d=1$, $A=B=C=0$,
\item[(ii)] $d\neq 1$, $A=0$, $B=a=(b-C)(d-1)$, $c=C/(d-1)$.
\end{enumerate}
\end{enumerate}

By \cite[Proposition 4.2]{CG3} every Ore-type twisted tensor product is a graded Ore extension of $\k[x,y]$. If $T$ is in normal form, $T\cong \k[x,y][z; \n, \d]$ where 
\begin{align*}
\n(x) &= dx+ey & \d(x)&=ax^2+bxy+cy^2\\
\n(y) &= Dy & \d(y) &= Ax^2+Bxy+Cy^2.
\end{align*}

\begin{thm}
\label{Ore types are semistandard}
An Ore-type twisted tensor product $T$ is semi-standard, and the following are equivalent.
\begin{enumerate}
\item $T$ is nondegenerate,
\item $T$ is AS-regular,
\item $\n$ is invertible.
\end{enumerate}
\end{thm}

To facilitate the proof, and for use below, we record the matrices $M$ and $N$ of equation (\ref{MandN}) of Section 2 and their determinants for an Ore-type twisted tensor product.

$$
M=\begin{bmatrix}
-y_0 & x_0 & 0\\
z_0-ax_0 & -bx_0-cy_0 & -dx_0-ey_0\\
-Ax_0 & z_0-Bx_0-Cy_0 &-Dy_0
\end{bmatrix}$$
$$
N=\begin{bmatrix}
y_1  & -ax_1-by_1-dz_1 & -Ax_1-By_1\\
-x_1 & -cy_1-ez_1 & -Cy_1-Dz_1\\
0  & x_1 & y_1
\end{bmatrix}
$$

\begin{align*}
\det M &=
dAx^3+(dB+eA-aD)x^2y+(eB+dC-bD)xy^2\\
&\quad +(eC-cD)y^3+(D-d)xyz-ey^2z\\
\det N &=
Ax^3+(B-a)x^2y+(C-b)xy^2-cy^3+(D-d)xyz-ey^2z
\end{align*}

\begin{proof}
The fact that $T$ is semi-standard, that is, $\det M= k\det N$ for some $k\in\k^*$, is easily verified in each case listed above.

The statement (3)$\Rightarrow$(2) follows from \cite[Theorem 6.2(i)]{CG3}, and (2)$\Rightarrow$(1) is \cite[Theorem 5.1]{ATVI}. Thus it suffices to prove (1)$\Rightarrow$(3).

If either $x$ or $y$ is nonzero, then $N$ clearly has rank 2. If $x=y=0$, then we may assume $z=1$ and  $N$ becomes
$$\begin{bmatrix}
0  & -d & 0\\
0 & -e & -D\\
0  & 0 & 0
\end{bmatrix}$$
which has rank 2 if and only if $dD\neq 0$. This condition implies the endomorphism $\n:\k[x,y]\to \k[x,y]$ is bijective in degree 1, so $\n$ is an isomorphism.
\end{proof}

The analysis of the cases O(1)(i,ii,iii,iv), O(2)(i,ii) is more involved than for the reducible-type twisted tensor products of Section \ref{reducible type ttps}. Consequently, following two preliminary results, this section is divided into four subsections. Our classification results are summarized in Table \ref{OreTable}, and in the following theorem. Double lines in Table \ref{OreTable} delineate the four subsections. 

\begin{table}[h]
\caption{Mapping of Ore-type TTP cases}
\label{OreTable}
\setlength\tabcolsep{2pt}
\begin{tabular}{|c|c|c|c|c|}
\hline
TTP Case & Subcase & Algebra & Condition & IM-type\\
\hline
\multirow{7}{*}{O(1)(i)}& \begin{tabular}{c}$A=0\neq B-a$\\ $(C-b)^2\neq 4c(a-B)$ \end{tabular} & $T(\a, \b, \g)^*$ & $\a+\b+\g\neq 0$ & $T_1$\\
\cline{2-5}
&\multirow{2}{*}{\begin{tabular}{c}$A=0\neq B-a$\\ $(C-b)^2=4c(a-B)$ \end{tabular}} & \multirow{2}{*}{$W(\g,\e)$} & $\e=0$ & $WL_2$\\
\cline{4-5}
&&& $\e=1$ & $WL_3$\\
\cline{2-5}
&\multirow{3}{*}{\begin{tabular}{c|c}\multirow{3}{*}{\begin{tabular}{c}$A=0$\\$a=B$\\$C-b=0\neq c$ \end{tabular}} & $a=b=0$\\ \cline{2-2} & $a\neq 0$\\ \cline{2-2} & $a=0\neq b$ \end{tabular}}  & $L(0,0)$& & $TL_1$\\
\cline{3-5}
&& $L(1,0)$& & $TL_2$\\
\cline{3-5}
&& $L(0,1)$& & $TL_4$\\
\cline{2-5}
&\multirow{2}{*}{\begin{tabular}{c}$A=0= B-a$\\ $C-b=0=c$ \end{tabular}} & $P(\e)$ & $\e=1$& $P_2$\\
\cline{3-5}
&& $P(\e)$ & $\e=0$ & \multirow{2}{*}{$P_1$}\\
\cline{1-4}
\noalign{\vskip\doublerulesep\vskip-\arrayrulewidth}
\cline{1-4}
\multirow{3}{*}{O(1)(iv)}&\multirow{3}{*}{$C=0$}& \multirow{3}{*}{$S(d,D)^{**}$} & $d=D\notin \{0,1\}$ & \\
\cline{4-5}
&&& $dD=0$ & not AS-regular\\
\cline{4-5}
&&& $d\neq D$, $d,D\notin\{0,1\}$ & \multirow{2}{*}{$S_1$}\\
\cline{1-4}
\noalign{\vskip\doublerulesep\vskip-\arrayrulewidth}
\cline{1-4}
\multirow{3}{*}{O(1)(ii)}&& \multirow{3}{*}{$S'(d,\e)$} & $d\neq 0=\e$ & \\
\cline{4-5}
&&&$d\neq 0$, $\e=1$ & $S'_1$\\
\cline{4-5}
&&& $d=0$ & not AS-regular\\
\hline
\hline
O(2)(i) & $a\neq 0$ &  $W$ && $T'$\\
\hline
\multirow{2}{*}{O(2)(ii)} && \multirow{2}{*}{ $W(d)$ }& $d\notin\{0,1\}$ & $WL_1$\\
 \cline{4-5}
 &&& $d=0$ & not AS-regular\\
\hline
\end{tabular}
*isomorphism classes of $T(\a,\b,\g)$ are parametrized by points $[\a:\b:\g]\in \P^2$ up to permutation of coordinates.

**isomorphism classes of $S(d,D)$ are parametrized by pairs $(d,D)$ up to permutation of coordinates.
\end{table}

\begin{thm}
\label{OreSummary}
Every Ore-type twisted tensor product is isomorphic as a graded algebra to one of the algebras listed in Table \ref{OreTable}. The algebras listed in Table \ref{OreTable} are AS-regular and are pairwise non-isomorphic, except as noted. The table gives the types of the AS-regular algebras. 
\end{thm}

The proof of Theorem \ref{OreSummary} can be found at the end of this section. 

\begin{rmk}
\label{complete Ore types}
Comparing the algebras listed in Table \ref{OreTable} with those of \cite[Theorem 3.1]{IM} shows that when ${\rm char}\ \k=0$, all algebras of types $T_1, WL_2, WL_3,$ and $TL_4$ are isomorphic to Ore-type twisted tensor products of $\k[x,y]$ and $\k[z]$.
\end{rmk}

We begin our classification of Ore-type twisted tensor products up to graded algebra isomorphism with case O(1). Recall that case O(1) has $e=0$. We start with a few simple observations.

\begin{lemma}\ 
\label{Ore case (1)(iii)}
\begin{enumerate}
\item Any twisted tensor product from case O(1)(iii)  is isomorphic as a graded algebra to one from case O(1)(ii).
\item Any twisted tensor product from case O(1)(iv) is isomorphic to one from case O(1)(iv) with $C=0$.
\end{enumerate}
\end{lemma}

\begin{proof}
Given a twisted tensor product from case O(1)(iii), the isomorphism that interchanges $x$ and $y$ and fixes $z$ yields an algebra from case O(1)(ii).

Given an algebra from case O(1)(iv) with $C=1$, the isomorphism that fixes $x$ and $y$ and maps $z\mapsto z+y/(1-D)$ yields an algebra from case O(1)(iv) with $C=0$.
\end{proof}

\subsection{Twisted tensor products from case O(1)(iv)}
\label{case (1) (iv)}
For $d, D \in \k$, $d\neq 1\neq D$,  define $$S(d,D)=\k\la x,y,z\ra/\la xy-yx, zx-dxz, zy-Dyz\ra.$$ 
The definition of $S(d,D)$ obviously makes sense for all values of $d$ and $D$. The restrictions are imposed to avoid redundancies with other families of algebras defined below, and to align with the values permitted in the corresponding case O(1)(iv).
\begin{prop}
\label{case (1)(iv)}
Any twisted tensor product from case O(1)(iv) is isomorphic as a graded algebra to an algebra from the family $S(d, D)$. Conversely, $S(d,D)$ is a twisted tensor product belonging to case O(1)(iv).
\end{prop}

\begin{proof}
Given a twisted tensor product $T$ from case O(1)(iv), by Lemma \ref{Ore case (1)(iii)}, we may assume $C=0$, so the relations of $T$ are
$$xy-yx, zx-ax^2-dxz, zy-Bxy-Dyz,$$
where $a(D-1)=B(d-1)$. In this case, $d\neq 1\neq D$. The  map that fixes $x$ and $y$, and sends $z$ to $\lambda x+z$ where $\lambda = \dfrac{a}{1-d}=\dfrac{B}{1-D}$ is an isomorphism from $T$ to $S(d,D)$. 
The converse is immediate from the description of case O(1)(iv).
\end{proof}

\begin{prop}
\label{S(d,D) point scheme}
The algebra $S(d,D)$ is AS-regular if and only if $dD\neq 0$. When $S(d,D)$ is AS-regular, the point scheme $E$ and associated automorphism $\s$ are as follows.
\begin{itemize}
\item[(1)] If $d=D$,  $E=\P^2$ and the automorphism $\s$ is diagonal. (type $P_1$)
\item[(2)] If $d\neq D$, $E$ is the union of three distinct lines, and the automorphism $\s$ stabilizes these lines. (type $S_1$)
\end{itemize}
In either case, the order of $\s$ is the least common multiple of the orders of $d, D\in \k^*$.

Moreover, the algebras $S(d,D)$ and $S(d',D')$ are graded isomorphic if and only if $(d',D') \in \{(d,D), (D,d)\}$.
\end{prop}

\begin{proof}
The statement that $S(d,D)$ is AS-regular if and only if $dD\neq 0$ follows from Proposition \ref{case (1)(iv)} and Theorem \ref{Ore types are semistandard}.
The matrix $M$ for $S(d,D)$ is
$$\begin{bmatrix}
-y_0  & x_0 & 0\\
z_0 & 0 & -dx_0\\
0  & z_0 & -Dy_0
\end{bmatrix}.$$
The point scheme of $S(d,D)$ is $E=\mathcal Z((D-d)xyz)$, which is reduced, and the automorphism $\s$ is given by
$$\s([x:y:z])=\begin{cases} [dx:dy:z] & x\neq 0\\ [0:Dy:z] & x=0\\ \end{cases}.$$  Statements (1), (2) and the statement regarding the order of $\s$ are now clear. 

The characterization of graded algebra isomorphism classes of the $S(d,D)$ is a special case of \cite[Example 4.10]{Mori}, so the proof is omitted.
\end{proof}

\subsection{Twisted tensor products from cases O(1)(ii), O(1)(iii)}
\label{case (1) (ii)}

By Lemma \ref{Ore case (1)(iii)} (1) there is no loss of generality in analyzing case O(1)(ii). For $\e, d \in \k$ such that $\e\in \{0,1\}$ and $d\neq 1$, we define: $$S'(d,\e)=\k\la x,y,z\ra/\la xy-yx, zy-yz, zx-dxz+\e y^2\ra.$$
As in the preceding subsection, we have imposed a restriction on $d$ to avoid redundancy and to align with the values permitted in the corresponding case.

\begin{prop}
\label{case (1)(ii)}
Any twisted tensor product from case O(1)(ii)  is isomorphic as a graded algebra to an algebra from the family $S'(d,\e)$. Conversely, $S'(d,\e)$ is a twisted tensor product belonging to case O(1)(ii).
\end{prop}

\begin{proof}
Given a twisted tensor product $T$ from case O(1)(ii), the relations of $T$ are
$$xy-yx, zy-yz, zx-ax^2-bxy-cy^2-dxz,$$
where $d\neq 1$. The graded algebra map which fixes $x$ and $y$ and sends $z$ to $\lambda x+\mu y+z$, where $\lambda=\dfrac{a}{1-d}$ and $\mu=\dfrac{b}{1-d}$, is an isomorphism from $T$ to the algebra with relations $xy-yx, zy-yz, zx-dxz+cy^2$. If $c\neq 0$, rescaling $y$ yields the algebra $S'(d,1)$.
The converse is immediate from the description of case O(1)(ii).
\end{proof}

\begin{prop}
\label{S_1 types}
An algebra $S'(d,\e)$ is AS-regular if and only if $d\neq 0$.
If $S'(d,\e)$ is AS-regular, then the point scheme $E$ and associated automorphism $\s$ are as follows.
\begin{itemize}
\item[(1)] If $\e=0$, $E$ is a union of three distinct lines, and $\s$ stabilizes the components. (type $S_1$)
\item[(2)] If $\e\neq 0$, $E$ is a union of a line and a conic meeting at two points, and $\s$ stabilizes the components and the intersection points. (type $S'_1$)
\end{itemize}
In either case, the order of $\s$ is the order of $d\in \k^*$.

 The algebras $S'(d,\e)$ and $S'(d',\e')$ are isomorphic if and only if $\e'=\e$ and $d'=d$ or, when $d\neq 0$, $d'=d^{-1}$. Moreover, $S'(d,0)$ is not isomorphic to an algebra belonging to the family $S(d', D')$.
\end{prop}

\begin{proof}
By definition of $S'(d,\e)$, the matrix $M$ is
$$\begin{bmatrix}
-y_0  & x_0 & 0\\
z_0 & \e y_0 & -dx_0\\
0  & z_0 & -y_0
\end{bmatrix}.$$
Proposition \ref{case (1)(ii)}
 and Theorem \ref{Ore types are semistandard} imply that $S'(d,\e)$ is AS-regular if and only if $d\neq 0$. In this case
the point scheme is $E=\mathcal Z((1-d)xyz+\e y^3)$, which is reduced since $d\neq 1$, and the automorphism $\s$ is given by
$$\s([x:y:z])=\begin{cases} [x:y:z] & y\neq 0 \\ [dx:0:z] & y= 0\\ \end{cases}.$$ The description of the point scheme of $S'(d,\e)$ and the automorphism $\s$ follows.

Next, we characterize $S'(d,\e)$ up to graded algebra isomorphism. For the ``if'' part of the statement, observe that if $d\neq 0$, then $S'(d,\e)\cong S'(d^{-1},\e)$ via the isomorphism given by switching $x$ and $z$ and rescaling $y\mapsto \sqrt{-d}y$.

Conversely, suppose $S'(d,\e)\cong S'(d', \e')$. Since $S'(d, \e)$ is AS-regular if and only if $d\neq 0$, we need only consider two cases: $dd'\neq 0$ and $d=d'=0$. First assume $dd'\neq 0$ so $S'(d,\e)$ and $S'(d',\e')$ are AS-regular of dimension 3, and thus are geometric algebras. By considering types, the first part of the Proposition implies $\e=\e'$. 

Let $\psi:\P^2\to \P^2$ be an automorphism such that $\psi\s=\s'\psi.$ Since $d\neq 1$, $\s$ and $\s'$ restrict to non-identity maps on $y=0$. As $\s$ and $\s'$ are the identity when $y\neq 0$,  $\psi$ stabilizes the $y=0$ component.

Let $[\psi_{ij}]\in GL_3$ be a matrix representing $\psi$, up to scaling. Since $\psi$ stabilizes the line $y=0$, we have $\psi_{21}=\psi_{23}=0$. Since $\s$ and $\s'$ fix $[1:0:0]$ and $[0:0:1]$, the condition $\psi\s=\s'\psi$ implies
$$[d'\psi_{11}:0:\psi_{31}]=[\psi_{11}:0:\psi_{31}]\quad\text{and}\quad [d'\psi_{13}:0:\psi_{33}]=[\psi_{13}:0:\psi_{33}].$$
Since $d'\neq 1$, we have $\psi_{11}\psi_{31}=\psi_{13}\psi_{33}=0$, and since $\psi$ is invertible, we have $\psi_{11}\psi_{33}\neq \psi_{13}\psi_{31}$. So either $\psi_{11}=\psi_{33}=0$ or $\psi_{13}=\psi_{31}=0$. In the former case, $\psi_{13}\neq0\neq\psi_{31}$, and evaluating $\psi\s=\s'\psi$ at $[1:0:1]$ yields
$$[\psi_{13}:0:d\psi_{31}]=[d'\psi_{13}:0:\psi_{31}],$$
whence $d'=d^{-1}$. A similar argument in the case $\psi_{13}=\psi_{31}=0$ yields $d=d'$. This completes the characterization of $S'(d,\e)$ in the case $dd'\neq 0$. 

For the case $d=d'=0$, note that the relation $zx=0$ holds in $S'(0,0)$. It is easy to check that products of linear elements in $S'(0,1)$ are always nonzero, hence $S'(0,\e)\cong S'(0,\e')$ implies $\e=\e'$.

Finally, observe that if the definition of $S(d,D)$ is extended to allow $D=1$, then $S'(d,0)=S(d,1)$. The fact that $S'(d,0)\ncong S(d',D')$ for $D'\neq 1$ is therefore a consequence of \cite[Example 4.10]{Mori}.
\end{proof}

\subsection{Twisted tensor products from case O(1)(i)}
\label{case (1) (i)}
Among the algebras from case O(1), it remains to consider case O(1)(i). Recall that in this case $d=D=1$; that is, the ring endomorphism $\n$ defining the Ore-extension $T$ is the identity. So, by Theorem \ref{Ore types are semistandard}, the algebra $T$ is AS-regular. As such, the classification in this case  is aided by consideration of the possible point schemes. The following reduction simplifies the point scheme calculations.

\begin{lemma}
\label{case (1)(i) reduction}
Any twisted tensor product from case O(1)(i) is isomorphic to one from case O(1)(i) with $A=0$.
\end{lemma}

\begin{proof}
Suppose $T$ is a twisted tensor product from case O(1)(i) with $A=1$. 
Let $p(t)=ct^3+(b-C)t^2+(a-B)t-1\in \k[t]$. 

If $p(t)$ is non-constant, let $\lambda\in \k$ be a zero of $p(t)$. Then the linear map $x\mapsto x$, $y\mapsto \l x+y$, $z\mapsto z$ determines an isomorphism of $T$ with $T'=\k\la x,y,z\ra/I'$ where $I'$ is generated by
$$xy-yx, zx-(a+b\l+c\l^2)x^2-(b+2c\l)xy-cy^2-xz,$$
$$zy-(B+(2C-b)\l-2c\l^2)xy-(C-c\l)y^2-yz.$$
Rescaling $y$, if necessary, this presentation satisfies the conditions of case O(1)(i).

If $p(t)$ is constant, then $c=b-C=a-B=0$. Interchanging $x$ and $y$ yields the isomorphic algebra with relations
$$xy-yx, zx-bx^2-axy-y^2-xz, zy-bxy-ay^2-yz.$$
Rescaling $y$ if necessary, this presentation satisfies the conditions of case O(1)(i).
\end{proof}

If $T$ is an algebra from case O(1)(i) with $A=0$, then the point scheme of $T$ is  
$$E=\mathcal Z(y\left[ (B-a)x^2+(C-b)xy-cy^2\right]).$$
The automorphism $\s$ is given by $$\s([x:y:z])=\begin{cases} [x:y:z+Bx+Cy] & y\neq 0\\ [x:0:ax+z] & y=0\\ \end{cases}.$$
We see that whenever $E$ consists of multiple components, the automorphism $\s$ stabilizes the components, and the components intersect in a single point $[0:0:1]$. 

When the quadratic form $$Q(x,y)=(B-a)x^2+(C-b)xy-cy^2$$ is a product of distinct linear factors, $E$ is either a union of three lines (when $a\neq B$), or the union of the double line $y^2=0$ and another line. We describe the three line case next.

For $\a, \b, \g \in \k$ such that $\a+\b+\g \ne 0$, define $$T(\a,\b,\g)=
\dfrac{\k\la x,y,z\ra}{\la xy-yx, xz-zx-\b x^2+(\b+\g)xy, yz-zy-\a y^2+(\a+\g)xy \ra}.$$

\begin{prop}
\label{T_1 types}
An algebra $T$ is a twisted tensor product from case O(1)(i) with $A=0$ and $B-a\neq 0\neq (C-b)^2-4c(a-B)$ if and only if $T$ is isomorphic to an algebra from the family $T(\a,\b,\g)$.
\end{prop}

\begin{proof}
Suppose $T$ is an algebra as hypothesized in the statement. Choose $\rho \in \k$ such that
$(a-B)^2+(4(a-B)c-(b-C)^2)\rho^2=0$ and let $\l=\dfrac{B-a+(b-C)\rho}{2(B-a)}$.
Then $x\mapsto x+\l y$, $y\mapsto \rho y$, $z\mapsto z$ determines an isomorphism from $T$ to an algebra $T(\a,\b,\g)$ where $\a=-(C\rho+B\l)$, $\b=-a$, and $\g =B(1+\l)+C\rho$.

Conversely, it is straightforward to check that when $\a+\b+\g\neq 0$, the algebra $T(\a,\b,\g)$ belongs to case O(1)(i) and satisfies the stated conditions.
\end{proof}

\begin{prop}
\label{point scheme of T_1}
The algebra $T(\a,\b,\g)$ is AS-regular. The point scheme of $T(\a,\b,\g)$ is a union of three distinct lines meeting at a point. The automorphism $\s$ stabilizes the components (type $T_1$), and the order of $\s$ is the least common multiple of the additive orders of $\a, \b, \g$ in $\k$. 

Moreover, $T(\a,\b,\g)\cong T(\a',\b',\g')$ if and only if the coordinates of the point $[\a': \b': \g']\in \P^2$ are a permutation of those of $[\a:\b:\g]$.
\end{prop}

\begin{proof}
As noted above, AS-regularity of $T(\a,\b,\g)$ follows from Theorem \ref{Ore types are semistandard}, since the automorphism $\n$ is the identity.
The matrix $M$ is
$$\begin{bmatrix}
-y_0 & x_0 & 0\\
z_0+\b x_0 & -(\b+\g)x_0 & -x_0\\
-(\a+\g)x_0 & z_0+\a x_0 &-y_0
\end{bmatrix},$$
so the point scheme of  is $E=\mathcal Z((\a+\b+\g)xy(x-y))$, which is reduced, and the automorphism is 
\begin{align*}
\s([0:y:z])&=[0:y:z+\a y]\\
\s([x:0:z])&= [x:0:z+\b x]\\
\s([x:x:z])&= [x:x:z-\g x].
\end{align*}
Noting that the lines $x=0$, $y=0$, and $x=y$ intersect in $[0:0:1]$, the first part of the Corollary follows.

Since the algebras $T(\a,\b,\g)$ are AS-regular of dimension 3, they are geometric algebras. Suppose $T(\a,\b,\g)\cong T(\a',\b',\g')$ and let $\s'$ be the automorphism of the point scheme $E'$ of $T(\a', \b',\g')$. Then there exists an automorphism $\psi:\P^2\to \P^2$ such that $\psi\s=\s'\psi$, and $\psi$ acts by a permutation on the triple of lines ($x=0$, $y=0$, $x=y$). Applying the inverse of this permutation to the triple $(\a, \b, \g)$ yields $(\a', \b', \g')$ up to a scalar. We verify this only in the case of a 3-cycle, the other cases follow by similar, and easier, arguments.

Assume that $\psi$ cyclically permutes the lines: $(x=0)\rightarrow (y=0)\rightarrow (x=y)\rightarrow (x=0)$. Let $[\psi_{ij}]\in GL_3$ represent the transformation $\psi$ up to rescaling. Without loss of generality we may assume 
$$[\psi_{ij}] = \begin{bmatrix} 1 & -1 & 0\\ 1 & 0 & 0 \\ \psi_{31} & \psi_{32} & \psi_{33}\\ \end{bmatrix}$$
where $\psi_{33}\neq 0$. Evaluating the commutativity relation $\psi\s=\s'\psi$ at $[0:1:0]$, $[1:0:0]$ and $[1:1:0]$ yields
$$\b'=-\a\psi_{33},\quad \g'=-\b\psi_{33},\quad \a'=-\g\psi_{33}.$$
Hence $[\a':\b':\g']=[\g:\a:\b]$ as desired.

Conversely, suppose that $[\a: \b: \g] \in \P^2$ and $\a+\b+\g \ne 0$. Let $E$, $E'$ denote the point schemes of $T(\a, \b, \g)$, $T(\g,\a,\b)$, with automorphisms $\s$, $\s'$, respectively. Define $\psi \in {\rm Aut}(\P^2)$ by the matrix $$[\psi_{ij}] = \begin{bmatrix} 1 & -1 & 0\\ 1 & 0 & 0 \\ 0 & 0 & -1\\ \end{bmatrix}.$$ Then one checks that $\psi \s = \s' \psi$, so, by Theorem \ref{geometric algebra isomorphism}, $T(\a, \b, \g) \cong T(\g, \a, \b)$. Similarly, $T(\a, \b, \g)$ is isomorphic to $T(\a', \b', \g')$ whenever $(\a', \b', \g')$ is any permutation of $(\a, \b, \g)$. 
\end{proof}
 
When $B-a=0\neq C-b$, the quadratic form $Q(x,y)$ has distinct linear factors, one of which is $y$. These algebras are isomorphic to those where the form is a perfect square and are considered below.

\begin{lemma}
\label{case (1)(i) reduction 2}
A twisted tensor product from case O(1)(i) with $A=0$ and $B-a=0\neq C-b$ is isomorphic to an algebra from case O(1)(i) with $A=0$ and $B-a\neq 0=(C-b)^2-4(a-B)c$.
\end{lemma}

\begin{proof}
An isomorphism is given by the mapping $x\mapsto y+\b x$, $y\mapsto x$, $z\mapsto z$, where $(C-b)\b-c=0$.
\end{proof}

When  $Q(x,y)$ is the nonzero perfect square $[2(B-a)x+(C-b)y]^2$, the point scheme of $T$ is either the union of this double line with $y=0$, or the triple line $y^3=0$. The next proposition simplifies the description of $T$ in the former case.

For $\g \in \k$ and $\e \in \{0,1\}$, define $$W(\g,\e)=
\dfrac{\k\la x,y,z\ra}{\left\la xy-yx, zx-xz+\e x^2+\g xy, zy-yz+\e xy+(1+\g)y^2\right\ra}.$$

\begin{prop}
\label{e=0 double lines}
An algebra $T$ is a twisted tensor product from case O(1)(i) with $A=0$ and $B-a\neq 0=(C-b)^2-4(a-B)c$ if and only if $T$ is isomorphic to an algebra from the family $W(\g, \e)$.
\end{prop}

\begin{proof}
Suppose $T$ is an algebra as hypothesized in the statement. The defining relations of $T$ are 
$$xy-yx, zx-ax^2-bxy-cy^2-xz, zy-Bxy-Cy^2-yz.$$
Let $\b$ be the unique solution to $(B-a)t^2+(C-b)t-c=0.$ It follows that $2(B-a)\b+(C-b)=0$.
The map sending $x\mapsto x+\b y$ and fixing $y$ and $z$ defines an isomorphism of $T$ to the algebra $T'$ with relations
$$xy-yx, zx-xz-ax^2-\l xy, zy-yz-Bxy-\l y^2,$$
where $\l=C+\b B$. 
If $\l \neq 0$, the map $x\mapsto -\dfrac{y}{a-B}$, $y\mapsto -\dfrac{x}{\l}$, $z \mapsto z$ defines an isomorphism from $T'$ to the algebra $W(\g, 1)$, where $\g = B/(a-B)$. 
If $\l =0$, then $x\mapsto -\dfrac{y}{a-B}$, $y\mapsto -x$, $z \mapsto z$ defines an isomorphism from $T'$ to the algebra $W(\g, 0)$, where $\g = B/(a-B)$. 

Conversely, for any $\g \in \k$ and $\e \in \{0,1\}$, $W(\g,\e)$ satisfies the conditions of case O(1)(i) with $A=0$, $B-a=0\neq C-b$. The result follows from Lemma \ref{case (1)(i) reduction 2}.
\end{proof}

We use the following result to characterize isomorphism classes of Ore-type algebras in the cases where the point scheme is not reduced.

Recall that if $\l \in \k^*$ and $\varphi: B \to B'$ is an isomorphism of graded $\k$-algebras, then the map $\varphi_{\l}: B \to B'$, given by $\varphi_{\l}(b) = \l^i b$ for all $b \in B_i$, is also an isomorphism of graded $\k$-algebras. We refer to $\varphi_{\l}$ as a scaling of $\varphi$.

\begin{prop}
\label{Ore isomorphism}
Let $A$ be a graded $\k$-algebra. Let $A[z, \n, \d]$, $A[w, \n', \d']$ be graded Ore extensions with $\deg z = \deg w = 1$. Let $\varphi:A[z;\n,\d]\to A[w;\n',\d']$ be a graded $\k$-algebra isomorphism.
\begin{itemize}
\item[(1)] If $\varphi(A)=A$, then there exists a scaling of $\varphi$, $\varphi_{\l}$, such that $\varphi_{\l}\n=\n'\varphi_{\l}$ and $\varphi_{\l}\d=\d'\varphi_{\l}$. 

\item[(2)] If $A=\k[x,y]$ and $Z(A[w;\n',\d'])_1=0$ then $\varphi(A)=A$.
\end{itemize}
\end{prop}

\begin{proof}
Assume that $\varphi(A)=A$. Since $\varphi$ is surjective and degree-0, there exists $k\in \k^*$ such that $\varphi(z)=a' + kw$ for $a'\in A_1$. Scaling $\varphi$, if necessary, we may assume that $k=1$. Then for any homogeneous $a\in A_i$ we have 
\begin{align*}
0 &=\varphi(za-\n(a)z-\d(a))\\
&=(a'+w)\varphi(a)-\varphi(\n(a))(w+a')-\varphi(\d(a))\\
&=[\n'(\varphi(a))-\varphi(\n(a))]w+a'\varphi(a)-\varphi(\n(a))a'+\d'(\varphi(a))-\varphi(\d(a)).
\end{align*}
Since $A[w;\n',\d']$ is a free left $A$-module on the basis $\{w^i\}$, $A$ is graded, $\deg \n= 0$, and $\deg \d = \deg \d' = 1$, the first statement follows. 

Now let $A=\k[x,y]$ and suppose that $Z(A[w;\n',\d'])_1=0$. 
Let $\n'(x)=m_{11}x+m_{12}y$ and $\n'(y)=m_{21}x+m_{22}y$.
Write $\varphi(x)=\a_1x+\b_1y+\g_1w$ and $\varphi(y)=\a_2x+\b_2y+\g_2w$. 
Then 
\begin{align*}
0=\varphi(xy-yx)&=[(\a_2\g_1-\a_1\g_2)(m_{11}-1)+(\b_2\g_1-\b_1\g_2)m_{21}] xw\\
&+[(\a_2\g_1-\a_1\g_2)(m_{12})+(\b_2\g_1-\b_1\g_2)(m_{22}-1)] yw\\
&+(\a_2\g_1-\a_1\g_2)\d'(x) + (\b_2\g_1-\b_1\g_2)\d'(y).
\end{align*}
Let $r = \a_2\g_1-\a_1\g_2$ and $s = \b_2\g_1-\b_1\g_2$. It follows from the calculation above and direct computation that $rx+sy$ is in $Z(A[w;\n',\d'])_1$, hence $r = s = 0$.
Thus both $\{(\a_1,\a_2),(\g_1,\g_2)\}$  and $\{(\b_1,\b_2),(\g_1,\g_2)\}$ are linearly dependent. If $\varphi(A) \not\subset A$, then $(\g_1,\g_2)\neq (0,0)$ and $\{\varphi(x),\varphi(y)\}$ is linearly dependent, a contradiction. 
\end{proof}

\begin{prop}
\label{double line algebras}
The algebras $W(\g,\e)$ and $W(\g',\e')$ are isomorphic if and only if $\e'=\e$ and $\g'=\g$.
\end{prop}

\begin{proof}
First, if $\e=1$, or $\e=0$ and $\g\notin\{0,-1\}$, then $\d(x) = -\e x^2-\g xy$ and $\d(y) = -\e xy-(1+\g)y^2$ are linearly independent, and so $Z(W(\g,\e))_1=0$. These $W(\g, \e)$ are not isomorphic to $W(0,0)$ or $W(-1,0)$ because $x$ or $y$ (respectively) is the unique central element in degree 1, up to scaling. Moreover, $W(0,0)$ and $W(-1,0)$ are not isomorphic since it is clear that $W(0,0)/\la x \ra$ is not isomorphic to $W(-1,0)/\la y \ra$. Henceforth we assume that $\e=1$, or $\e=0$ and $\g\notin\{0,-1\}$.

Suppose there is an isomorphism $\varphi:W(\g,\e)\to W(\g',\e')$. Then $\d(x)$ and $\d(y)$ are linearly independent elements of $W(\g,\e)$ as are their counterparts $\d'(x)$ and $\d'(y)$ in $W(\g',\e')$. By Proposition \ref{Ore isomorphism} we may write $\varphi(x)=\a_1x+\b_1y$ and $\varphi(y)=\a_2x+\b_2y$, and we may assume that $\varphi$ commutes with $\d$ and $\d'$. Equating coefficients of $x^2$, $xy$ and $y^2$ in $\varphi(\d(x))=\d'(\varphi(x))$ and $\varphi(\d(y))=\d'(\varphi(y))$ yields the following:
\begin{align*}
\a_1(\e\a_1+\g\a_2) &= \a_1\e' & \a_2(\e\a_1+(1+\g)\a_2) &= \a_2\e'\\ 
\b_1(\e\b_1+\g\b_2) &= \b_1(1+\g') & \b_2(\e\b_1+(1+\g)\b_2) &= \b_2(1+\g')
\end{align*}
$$2\e\a_1\b_1+\g(\a_1\b_2+\a_2\b_1) =\a_1\g'+\b_1\e' $$ 
$$2(1+\g)\a_2\b_2+\e(\a_1\b_2+\a_2\b_1) =\a_2\g'+\b_2\e'.$$

Suppose that $\e=0$ and $\e'=1$. Considering the first two equations, we see that $\a_1=0$. The fifth and second equations then imply that $\b_1=0$,  so $\varphi(x)=0$, a contradiction. Thus we must have $\e=\e'$. 

When $\e=\e'=0$, then $\a_2=0$ by the second equation, so $\a_1\neq 0\neq \b_2$, else $\varphi$ is not surjective. The fourth and fifth equations jointly imply that $\b_2=1$, so $\g=\g'$.

When $\e=\e'=1$, we also have $\a_2=0$. If $\a_1\neq 0$, this follows from the first pair of equations. If 
 $\a_1=0$, then $\b_1\neq 0$ and $(1+\g)\a_2=1$. The fifth equation then simplifies to $\g\a_2=1$, which implies that $\a_2 = 0$. A straightforward check shows that $\a_2=0$ implies $\g=\g'$.
\end{proof}

\begin{prop}
\label{double line point scheme}
The algebra $W(\g,\e)$ is AS-regular. The point scheme of $W(\g,\e)$ is the union of a line and a double line. The automorphism $\s$ stabilizes the components, and 
\begin{itemize}
\item[(1)] if $\e=0$, then $\s$ restricts to the identity on the double line (type $WL_2$); 
\item[(2)] if $\e=1$, then $\s$ has order ${\rm char}\ \k$ on the double line. (type $WL_3$)
\end{itemize}
\end{prop}

\begin{proof} The algebra $W(\g,\e)$ is AS-regular by Proposition \ref{e=0 double lines} and Theorem \ref{Ore types are semistandard}. The matrix $M$ is
$$
\begin{bmatrix}
-y_0 & x_0 & 0\\
z_0+\e x_0 & \g x_0 & -x_0\\
0 & z_0+\e x_0+(1+\g)y_0 &-y_0
\end{bmatrix}$$
so the point scheme is $E=\mathcal Z(xy^2)$ and the automorphism is 
\begin{align*}
\s([0:y:z])&=[0:y:z+(1+\g)y]\\
\s([x:0:z])&=[x:0:z+\e x].
\end{align*}
The result follows immediately.
\end{proof}

For $\e_1, \e_2\in \{0,1\}$ we define
$$L(\e_1,\e_2)=\dfrac{\k\la x,y,z\ra}{\la xy-yx, zx-xz+\e_1xy+\e_2x^2, zy-yz+x^2+\e_1y^2+\e_2xy\ra}.$$

\begin{prop}
\label{triple lines}
An algebra $T$ is a twisted tensor product from case O(1)(i) with $A=0$ and $B-a=C-b=0\neq c$ if and only if $T$ is isomorphic to one of $L(0,0), L(0,1),$ or $L(1,0)$, no two of which are isomorphic.
\end{prop}

\begin{proof}
Suppose $T$ is an algebra as hypothesized in the statement. If $a\neq 0$, then $x\mapsto y-(b/\sqrt{ac})x$, $y\mapsto \sqrt{a/c}x$, $z\mapsto -az$
determines an isomorphism from $T$ to $L(1,0)$.
If $a=b=0$, then $x\mapsto -cy$, $y\mapsto x$, $z\mapsto z$ determines an isomorphism from $T$ to $L(0,0)$.
If $a=0\neq b$, then $x\mapsto (c/b)y$, $y\mapsto x$, $z\mapsto -bz$ determines an isomorphism from $T$ to $L(0,1)$. Conversely, it is clear that $L(0,0)$, $L(0,1)$ and $L(1,0)$ are algebras from case O(1)(i) of the stated form.

To see that the three algebras are pairwise non-isomorphic, first note that the degree-1 generator  $x\in L(0,0)$ is central, and $L(1,0)$ and $L(0,1)$ contain no central elements in degree 1. By Proposition \ref{Ore isomorphism}, any isomorphism $\varphi:L(1,0)\to L(0,1)$ must restrict to an automorphism of $\k[x,y]$ and commute with the respective derivations. An easy direct calculation in degree 1 shows that no such $\varphi$ exists.
\end{proof}

\begin{prop}
\label{triple line point scheme}
The algebras $L(0,0)$, $L(0,1)$, and $L(1,0)$ are AS-regular, and the point scheme of each is a triple line. The order of the automorphism $\s$ is the order of $\e_1$ in the additive group $(\k, +)$. 
\end{prop}

\begin{proof} The algebras $L(0,0)$, $L(0,1)$, and $L(1,0)$ are AS-regular by Proposition \ref{triple lines} and Theorem \ref{Ore types are semistandard}. The matrix $M$ is
$$
\begin{bmatrix}
-y_0 & x_0 & 0\\
z_0+\e_2x_0 & \e_1x_0 & -x_0\\
x_0 & z_0+\e_2x_0+\e_1y_0 &-y_0
\end{bmatrix}$$
so the point scheme of $T$ is $E=\mathcal Z(x^3)$ and the automorphism $\s$ is given by 
$\s([0:y:z])=[0:y:z+\e_1y].$ The result follows.
\end{proof}

\begin{rmk}
\label{triple line types}
The presentations of $L(0,0)$, $L(1,0)$, and $L(0,1)$ align with those of Types $TL_1, TL_2,$ and $TL_4$ given in \cite[Theorem 3.1]{IM}., respectively. We remark that $\s$ is the identity on closed points for both $L(0,0)$ and $L(0,1)$.
\end{rmk}

For $\e \in \{0, 1\}$ define $$P(\e)=\k\la x,y,z\ra/\la xy-yx, zx-\e x^2-xz, zy-\e xy-yz\ra.$$
\begin{prop}
\label{full P^2 cases}
An algebra $T$ is a twisted tensor product from case O(1)(i) with $B-a=C-b=c=0$ if and only if $T$ is isomorphic to $P(\e)$. 
\end{prop}

\begin{proof}
Suppose $T$ is an algebra as hypothesized in the statement. The defining relations of $T$ are
$$xy-yx, zx-ax^2-bxy-xz, zy-axy-by^2-yz.$$
If $a=b=0$, then clearly $T=P(0)$. 
If $a\neq 0$, then the change of variables $x\mapsto x-(b/a)y$, $y \mapsto y$, $z\mapsto az$  is an isomorphism from $T$ to the algebra $P(1)$. 
If $a=0\neq b$, then interchanging $x$ and $y$ yields an algebra of the form just considered. 

Conversely, the algebras $P(0)$ and $P(1)$ clearly satisfy the conditions specified in this subcase of case O(1)(i).
\end{proof}

\begin{prop}
\label{point scheme is P2}
The algebras $P(0)$ and $P(1)$ are AS-regular. The point scheme of each is $\P^2$. The automorphism $\s$ is the identity for $P(0)$, and for $P(1)$ the order of $\s$ is the characteristic of $\k$.
\end{prop}

\begin{proof} The algebras $P(0)$ and $P(1)$ are AS-regular by Proposition \ref{full P^2 cases} and Theorem \ref{Ore types are semistandard}. The matrix $M$ is
$$
\begin{bmatrix}
-y_0 & x_0 & 0\\
z_0-\e x_0 & 0 & -x_0\\
0 & z_0-\e x_0 &-y_0
\end{bmatrix}$$
so the point scheme of $P(\e)$ is $\mathcal Z(0)=\P^2$ and the automorphism $\s$ is given by 
$\s([x:y:z])= [x:y:z-\e x]$. 
\end{proof}

\subsection{Twisted tensor products from cases O(2)(i), O(2)(ii)}
\label{case (2}
Now we turn our attention to algebras from case O(2). Recall that in this case $e=1$. Define
$$W=\k\la x,y,z\ra/\la xy-yx, zx-xz+x^2+yz, zy-yz\ra.$$

\begin{prop}
\label{case (2)(i)}
Any twisted tensor product from case O(2)(i) with $a=0$ is isomorphic to an algebra from case O(1)(i). 

Any twisted tensor product from case O(2)(i) with $a\neq 0$ is isomorphic to  $W$.
\end{prop}

\begin{proof}
Let $T$ be a twisted tensor product from case O(2)(i). The defining relations of $T$ are
$$xy-yx, zx-ax^2-bxy-cy^2-xz-yz, zy-yz.$$
If $a=0$, the map determined by $x\mapsto -z$, $y\mapsto y$, and $z\mapsto x-cy+bz$ is an isomorphism from $T$ to the algebra with relations
$$xy-yx, zx-xy-xz, zy-yz,$$
which belongs to case O(1)(i).

If $a\neq 0$, the isomorphism given by $x\mapsto -x$, $y\mapsto y$, $z\mapsto az-cy+bx$ identifies $T$ with the algebra $W$ defined above.
\end{proof}

\begin{prop}
\label{point scheme of W}
The algebra $W$ is AS-regular. The point scheme of $W$ is the union of a line and an irreducible conic. The components meet at a single point, and $\s$ stabilizes the components (type $T'$). The order of $\s$ is ${\rm char}\ \k$.
\end{prop}

\begin{proof} By Proposition \ref{case (2)(i)} and Theorem \ref{Ore types are semistandard}, $W$ is AS-regular.
The matrix $M$ is
$$
\begin{bmatrix}
-y_0 & x_0 & 0\\
z_0+x_0 & 0 & -x_0+y_0\\
0 & z_0 &-y_0
\end{bmatrix}$$
so the point scheme of $W$ is $\mathcal Z(y(x^2+yz))$. The components intersect in the point $[0:0:1]$, and the automorphism is given by
$$\s([x:y:z])=\begin{cases} [x:y:z] & y\neq 0\\ [x:0:x+z] & y=0\\ \end{cases}.$$ The result follows.
\end{proof}

For $d\in \k$, $d\neq 1$, 
define
$$W(d)=\k\la x,y,z\ra/\la xy-yx, zx-dxz-yz, zy-dyz\ra.$$

\begin{prop}
\label{case (2)(ii)}
A twisted tensor product $T$ belongs to case O(2)(ii) if and only if $T$ is isomorphic to $W(d)$ for some $d\neq 1$. 
Moreover, the algebras $W(d)$ and $W(d')$ are isomorphic if and only if $d'=d$.
\end{prop}

\begin{proof}
If $T$ is a twisted tensor product from case O(2)(ii), then the map $x \mapsto x$, $y \mapsto y$,  $z \mapsto z+\l x-cy$, where $\l=\dfrac{a}{1-d}$, is an isomorphism from $T$ to an algebra in the family $W(d)$. Conversely, it is clear that $W(d)$ is an algebra from case O(2)(ii).

For the last statement, since $d\neq 1$, one checks that $Z(W(d))_1=0$, so by Proposition \ref{Ore isomorphism}, an isomorphism $\varphi:W(d)\rightarrow W(d')$ would restrict to an automorphism of $\k[x,y]$ that commutes with $\n$ and $\n'$. As the Jordan forms of $\n_1$ and $\n'_1$ are indecomposable, the condition $\varphi\n=\n'\varphi$ implies that $d=d'$. 
\end{proof}

\begin{prop}
\label{point scheme of W(d)}
The algebra $W(d)$ is AS-regular if and only if $d\neq 0.$ If $W(d)$ is AS-regular, then the point scheme is the union of a line and a double line. The automorphism $\s$ restricts to the identity on the reduced line. When restricted to the double line, the order of $\s$ is the additive order of $d\in\k$. (type $WL_1$)
\end{prop}

\begin{proof} By Proposition \ref{case (2)(ii)} and Theorem \ref{Ore types are semistandard}, $W(d)$ is AS-regular if and only if $d\neq 0$. The matrix $M$ is
$$
\begin{bmatrix}
-y_0 & x_0 & 0\\
z_0 & 0 & -dx_0-y_0\\
0 & z_0 &-dy_0
\end{bmatrix}$$
so the point scheme of $W(d)$ is $\mathcal Z(y^2z)$. The automorphism $\s$ is given by
\begin{align*}
\s([x:0:z]) &=[dx:0:z] \\
\s([x:y:0]) &=[x:y:0],
\end{align*}
and the result follows.
\end{proof}

We conclude this section with a proof of Theorem \ref{OreSummary}.

\begin{proof}[Proof of Theorem \ref{OreSummary}]
First we consider case O(1)(i). In this case, every algebra is AS-regular by Theorem \ref{Ore types are semistandard}. Lemma \ref{case (1)(i) reduction} shows there is no loss of generality in assuming $A=0$. We further divide case O(1)(i) into two subcases. Algebras from case O(1)(i) with $A=0\neq B-a$ are classified up to isomorphism as $T(\a,\b,\g)$ or $W(\g,\e)$ in Propositions \ref{T_1 types} and \ref{e=0 double lines}. Proposition \ref{point scheme of T_1} shows that $T(\a,\b,\g)$ is of Type $T_1$, and Proposition \ref{double line point scheme} shows that $W(\g,\e)$ is of Type $WL_2$ when $\e=0$, and Type $WL_3$ when $\e=1$.

Continuing with case O(1)(i), when $A=0=B-a$, Lemma \ref{case (1)(i) reduction 2} shows that no generality is lost in assuming $b=C$. Case O(1)(i) algebras with $A=B-a=C-b=0$ are classified up to isomorphism as $L(0,0), L(0,1), L(1,0)$, or $P(\e)$ in Propositions \ref{triple lines} and \ref{full P^2 cases}. Proposition \ref{triple line point scheme} and Remark \ref{triple line types} show that $L(0,0), L(1,0),$ and $L(0,1)$ are AS-regular of Type $TL_1$, $TL_2$, and $TL_4$, respectively. Proposition \ref{point scheme is P2} proves that $P(\e)$ is AS-regular of Type $P_1$ when $\e=0$, and Type $P_2$ when $\e=1$.

Algebras from case O(1)(ii) are isomorphic to algebras from the family $S'(d,\e)$ by Proposition \ref{case (1)(ii)}. Proposition \ref{S_1 types} shows $S'(d,\e)$ is AS-regular if and only if $d\neq 0$, and that the AS-regular algebras in this family are of Type $S_1$ when $\e=0$, and Type $S'_1$ when $\e=1$.

By Lemma \ref{Ore case (1)(iii)}(1), every algebra from case O(1)(iii) is isomorphic to an algebra from case O(1)(ii). 

For case O(1)(iv), Lemma \ref{Ore case (1)(iii)}(2) shows that no generality is lost in assuming $C=0$. These algebras are isomorphic to algebras from the family $S(d,D)$ by Proposition \ref{case (1)(iv)}. Proposition \ref{S(d,D) point scheme} shows $S(d,D)$ is AS-regular if and only if $dD\neq 0$, and proves that the AS-regular algebras in this family are of Type $P_1$ when $d=D$, and of Type $S_1$ when $d\neq D$. Proposition \ref{point scheme is P2} shows that the Type $P_1$ algebras isomorphic to $S(d,D)$ are not isomorphic to the algebra $P(0)$. Proposition \ref{S_1 types} shows that the Type $S_1$ algebras isomorphic to $S(d,D)$ are not isomorphic to Type $S_1$ algebras from the family $S'(d,0)$.

In case O(2)(i), Proposition \ref{case (2)(i)} shows that no generality is lost in assuming that $a\neq 0$ and shows that all algebras in this family are isomorphic to the algebra $W$. Proposition \ref{point scheme of W} establishes that $W$ is AS-regular of Type $T'$.

Finally, algebras from case O(2)(ii) are isomorphic to algebras from the family $W(d)$ by Proposition \ref{case (2)(ii)}. Proposition \ref{point scheme of W(d)} shows that $W(d)$ is AS-regular if and only if $d\neq 0$, and that the regular algebras in this family are of Type $WL_1$.
\end{proof}

\begin{rmk}\ 
\label{complete Ore cases}
\cite[Theorem 3.2]{IM} characterizes quadratic AS-regular algebras of non-EC type up to Morita equivalence. Types $P, T, T', WL, TL$, each of which include algebras from Theorem \ref{OreSummary}, consist of a single Morita equivalence class. For Type $S$, \cite[Theorem 3.2]{IM} implies that every Morita equivalence class contains an algebra with $\b=\g=1$, $\a\neq 0,1$. Such an algebra is the twisted tensor product $S'(\a,0)$ (see Subsection \ref{case (1) (ii)}). Similarly for Type $S'$, every Morita class contains an algebra with $\b=1$ and $\a\neq 0,1$. After interchanging $x$ and $y$ and rescaling $y$, we obtain the twisted tensor product $S'(1/\a,1)$.
\end{rmk}

\section{Graded Twisted Tensor Products and Sklyanin Algebras}
\label{Sklyanins}

As the preceding sections show, the only graded twisted tensor products of $\k[x,y]$ and $\k[z]$ whose point schemes are elliptic curves are the algebras $T(g,h)$ ($h\neq 0$) of Section \ref{elliptic type ttps}. These algebras belong to subtype B in the classification schemes of \cite{ATVI} and \cite{Mat}. It is therefore natural to ask whether other subtypes of type-EC AS-regular algebras are graded twisted tensor products of the form $A\tsr \k[z]$, where $A$ is a two-dimensional quadratic AS-regular algebra. We answer this question for the well-known family of three-dimensional Sklyanin algebras. 

By definition, a \emph{three-dimensional Sklyanin algebra} is any algebra of the form
$$S(a,b,c) = \dfrac{\k\la x,y,z\ra}{\la ayz+bzy+cx^2, azx+bxz+cy^2, axy+byx+cz^2\ra}$$
such that $[a:b:c]\in \P^2 - \mathcal{D}$, where $$\mathcal{D} = \{[1:0:0], [0:1:0], [0:0:1]\} \cup \{[a:b:c] : a^3 = b^3 = c^3\}.$$ This family is defined over a field $\k$ of any characteristic.  It is well known that $S(a,b,c)$ is AS-regular if and only if $[a:b:c] \notin \mathcal{D}$. If $[a: b: c] \in \mathcal{D}$, the algebra $S(a,b,c)$ is referred to as a \emph{degenerate Sklyanin algebra}. Let $$S_1  = \dfrac{\k\la u, v, w \ra}{\la u^2, v^2, w^2 \ra}, \qquad S_2 =   \frac{\k\la u, v, w \ra}{\la uv, vw, wu \ra}.$$ Smith \cite{Smith} proved that if $[a:b:c] \in \mathcal{D}$ and ${\rm char \ } \k \neq 3$, then $S(a, b, c) \cong S_1$ if $a = b$ and $S(a, b, c) \cong S_2$ if $a \ne b$.

The point scheme of the Sklyanin algebra $S(a,b,c)$ is 
$$E= \mathcal Z((a^3+b^3+c^3)xyz-(abc)(x^3+y^3+z^3)),$$
which describes an elliptic curve if and only if $abc\neq 0$ and $(3abc)^3\neq (a^3+b^3+c^3)^3$. When $E$ is an elliptic curve, the automorphism of $E$ can be described as translation by $[a:b:c]$ in the group law, where the identity element is $O_E = [1:-1:0]$ (see \cite[p. 38]{ATVI}).
In this section we are concerned mainly with those $S(a,b,c)$ whose point scheme is an elliptic curve. We say that $S(a,b,c)$ is of type EC if this is the case.

If $S(a,b,c)\cong A\tsr \k[z]$, then $S(a,b,c)$ contains a subalgebra isomorphic to $A$, generated in degree 1. We begin by characterizing which $S(a,b,c)$ of type EC contain a one-generated subalgebra isomorphic to a skew polynomial ring $\k_q[r,s]$, $q \in \k^*$, or the Jordan plane $\k\la r,s \ra/\la rs-sr+s^2\ra$.  

\subsection{Three-dimensional Sklyanin algebras containing a quantum $\P^1$}

If $A=T(V)/\la R\ra$ is a quadratic algebra, the \emph{quadratic dual} algebra is defined to be the algebra $A^!=T(V^*)/\la R^{\perp}\ra$ where $R^{\perp}\subset V^*\tsr V^*$ is the orthogonal complement to $R$ with respect to the natural pairing. Note that $(A^!)^!$ is canonically isomorphic to $A$. We establish a few facts about the schemes $\G(S(a,b,c)^!)$. Let $e_1 = [1:0:0]$, $e_2 = [0:1:0]$, $e_3 = [0:0:1]$.

\begin{prop}
\label{dualPointScheme}
Let $A = S(a,b,c)$ be a three-dimensional Sklyanin algebra.
\begin{itemize}
\item[(1)]  If $c \neq 0$, then $\G(A^!)$ is the empty scheme.
\item[(2)] $\G(S(a, -a, 0)^!)$ is the empty scheme. 
\item[(3)] $\G(S_1^!) = \{e_1 \times e_1, e_2 \times e_2, e_3 \times e_3\}.$
\end{itemize}
\end{prop}

\begin{proof} (1) Assume that $c \neq 0$. Since $\G(A^!)$ is a projective scheme it is quasi-compact; as such, to prove the result, it suffices to show that $\G(A^!)$ has no closed points. Since $c \ne 0$, we can write the defining relations of $A^!$ as
$$ax^2-cyz\qquad  ay^2-czx\qquad  az^2-cxy\qquad
bx^2-czy\qquad  by^2-cxz\qquad  bz^2-cyx.$$
Multilinearizing these relations, we see that a closed point 
$$[x_0:y_0:z_0]\times [x_1:y_1:z_1]\in \G(A^!)$$ satisfies
\begin{align*}
ax_0x_1-cy_0z_1 &= 0 \quad (E_1) & bx_0x_1-cz_0y_1 &= 0 \quad (E_2)\\ 
ay_0y_1-cz_0x_1 &= 0 \quad (E_3) & by_0y_1-cx_0z_1 &= 0 \quad (E_4)\\ 
az_0z_1-cx_0y_1 &= 0 \quad (E_5)& bz_0z_1-cy_0x_1 &= 0 \quad (E_6).
\end{align*}

First assume $a = 0$. Since $[a:b:c] \notin \mathcal{D}$, we have $bc \ne 0$. It is easy to check that no point in $\P^2\times \P^2$ satisfies all of the equations $(E_i)$. The case $b=0$ is analogous.

Now assume that $ab \neq 0$. Let $s: \P^2 \times \P^2 \to \P^8$, $$s([x_0: y_0: z_0],[x_1: y_1: z_1]) = \begin{bmatrix} x_0x_1 & x_0y_1 & x_0z_1 \\ y_0x_1 &y_0y_1 & y_0z_1 \\ z_0x_1 & z_0y_1 & z_0z_1 \end{bmatrix}$$ be the Segre embedding. Let $w_{ij}$, $1 \leq i, j \leq 3$ be coordinates on $\P^8$. Recall that the image $s(\P^2 \times \P^2)$ is the subscheme of $\P^8$ cut out by the $2 \times 2$ minors of the matrix $[w_{ij}]$. Equations $(E_i)$, $1 \leq i \leq 6$ show that the closed points of $s(\G(A^!))$ are the rank-1 matrices of the form $$M = \begin{bmatrix} c w_{11} & aw_{33} & bw_{22} \\   b w_{33} & cw_{22} & aw_{11} \\  a w_{22} & bw_{11} & cw_{33} \end{bmatrix}.$$ Let $M_{ij}$, $1 \leq i, j \leq 3$, denote the minors of $M$. If $M \in s(\G(A^!))$, then
consideration of the equations $w_{ii}M_{ii}=0$, for $i=1,2,3$
shows that $w_{11}^3 = w_{22}^3 = w_{33}^3$. Note that since ${\rm rank}(M) = 1$, $w_{11}w_{22}w_{33} \neq 0$. Furthermore, consideration of $cw_{33}M_{33}=0$, $bw_{11}M_{32}=0$, and $aw_{22}M_{31}=0$
shows that $$a^3 w_{11}w_{22}w_{33} = b^3 w_{11}w_{22}w_{33} = c^3 w_{11}w_{22}w_{33}.$$ Since $w_{11}w_{22}w_{33} \neq 0$ we conclude that $a^3 = b^3 = c^3$. This contradicts the assumption that $S(a,b,c)$ is a three-dimensional Sklyanin algebra.

For (2), the algebra $S(a, -a, 0)$ is the polynomial algebra on three variables. Its quadratic dual is an exterior algebra, and it is straightforward to check that there are no closed points in $\G(S(a, -a, 0)^!)$. 

For (3), a closed point $[x_0: y_0: z_0] \times [x_1:y_1:z_1] \in \G(S_1^!)$ satisfies $$x_0y_1 = y_0x_1 = z_0x_1 = x_0z_1 = z_0y_1 = y_0z_1 =0.$$ It is then straightforward to verify that  $\G(S_1^!) = \{e_1 \times e_1, e_2 \times e_2, e_3 \times e_3\}.$
\end{proof}

The following theorem characterizes the $S(a,b,c)$ of type EC that contain a one-generated skew polynomial ring.

\begin{thm}
\label{not a skew polynomial ring}
Let $S(a,b,c)$ be a three-dimensional Sklyanin algebra of type EC. Assume that ${\rm char}\ \k\neq 2,3$. 
\begin{itemize} 
\item[(1)] If $r, s\in S(a,b,c)_1$ are linearly independent and $rs=qsr$ for $q\in \k^*$, then $q=-1$ and $a=b$. 
\item[(2)] If the algebra $S(1,1,c)$ is of type EC, then there exist $r, s \in S(1,1,c)_1$ that are linearly independent and $rs = -sr$. More precisely, up to scaling in each component independently, the set of such pairs $(r, s)$ is in one-to-one correspondence with the closed points of the scheme $\G(S(-c,-c,2))$. 
\end{itemize}
\end{thm}

\begin{proof}
(1) Suppose that $r, s\in S(a,b,c)_1$ are linearly independent and $rs=qsr$ for some $q\in \k^*$. Since $abc\neq 0$, there is no loss of generality in assuming $a=1$. We write $r=\a_1x+\b_1y+\g_1z$ and $s=\a_2x+\b_2y+\g_2z$ for $\a_i, \b_i, \g_i\in \k$. The equation $rs=qsr$ implies the following identities hold in $\k$:
\begin{align*}
& (1-q)\a_1\a_2-c\b_1\g_2+qc\g_1\b_2=0 \quad  (E_1) \\ 
& (1-q)\b_1\b_2-c\g_1\a_2+qc\a_1\g_2=0 \quad  (E_2)\\ 
& (1-q)\g_1\g_2-c\a_1\b_2+qc\b_1\a_2=0 \quad  (E_3)\\ 
& (1+bq)\b_1\a_2-(q+b)\a_1\b_2=0 \quad  (E_4)\\
& (1+bq)\a_1\g_2-(q+b)\g_1\a_2=0\quad  (E_5)\\
& (1+bq)\g_1\b_2-(q+b)\b_1\g_2=0 \quad  (E_6).
\end{align*}

We refer to this set of equations as \emph{system $(E)$}. We view system $(E)$ as the defining equations of $\G(A)$, where $A$ is the quotient of  $\k\la x, y, z \ra$ by the relations:
\begin{align*}
&(1-q)x^2 - cyz+qczy, \qquad (1-q)y^2 - czx+qcxz, \qquad (1-q)z^2 - cxy+qcyx, \\
&(1+bq)yx-(q+b)xy, \qquad (1+bq)xz-(q+b)zx, \qquad (1+bq)zy-(q+b)yz.
\end{align*} Then we identify $r$ and $s$ with components of a closed point $p\in\G(A)$. Observe that under this identification, if $p\in \Delta(\P^2)$, where $\Delta :  \P^2 \to \P^2 \times \P^2$ is the diagonal embedding, then $r$ and $s$ are linearly dependent. We consider several cases.

First, if $q=1$, then $A = \k[x,y,z]$ is the polynomial ring on three variables. Since $\G(\k[x,y,z])=\Delta(\P^2)$, any $r$ and $s$ whose coefficients satisfy system $(E)$ must be linearly dependent, a contradiction. 
We assume henceforth that $q\neq 1$.

If $1+bq = q+b = 0$, then $q = -1$ and $b = 1$, and we have the conclusion of (1). So we assume, henceforth, that $1+bq \neq 0$, or $b+q \neq 0$. Under this assumption it is clear that the six defining relations of $A$ are linearly independent. Consider the quadratic dual $A^!$; one checks that $A^! = S(1+bq, q+b, c(1+q)).$
It is possible that this algebra is a degenerate Sklyanin algebra. 

If $q = -1$ and $b \neq 1$, then Proposition \ref{dualPointScheme}(2) shows that system $(E)$ has no solutions, a contradiction. Hence, for the remainder of the proof, we assume that $q \neq -1$. If we suppose that $A^!$ is not a degenerate Sklyanin algebra, then Proposition \ref{dualPointScheme}(1) shows that system $(E)$ has no solutions, a contradiction.  

Finally, we rule out the possibility that $A^!$ is a degenerate Sklyanin algebra. Suppose that this is the case. Note that $c(1+q) \neq 0$, and recall that $1+bq \neq 0$, or $b+q \neq 0$, so $[1+bq: b+q: c(1+q)] \notin \{e_1, e_2, e_3\}$. Hence, without loss of generality, $$(1+bq)^3 = (b+q)^3 = (c(1+q))^3 = 1.$$ If $1+bq = b+q$, then $A^! \cong S_1$, and Proposition \ref{dualPointScheme}(3) implies that $\G(A)$ is contained in $\Delta(\P^2)$, so $r$ and $s$ are linearly dependent, a contradiction. If $1+bq \ne b+q$, then $A^! \cong S_2$. Let $\a = 1+bq$, $\b = b+q$, $\g = c(1+q)$. Define elements of $S(1, b,c)$ by $$r = \a^{-1}x+\b^{-1}y+\g^{-1}z, \qquad s = \b^{-1}x+\a^{-1}y+\g^{-1}z, \qquad t = \a\b\g x+ \a\b\g y+ z.$$ It is easy to check that $\{r, s, t\}$ is linearly independent and, using system $(E)$, $$rs - qsr = rt-qrt = st - qts = 0.$$ Hence $S(1,b,c) \cong S(1, -q, 0)$. However, the point scheme of the algebra $S(1, -q, 0)$ is not an elliptic curve, a contradiction. 

For (2), we assume $S(1,1,c)$ is a Sklyanin algebra of type EC. Therefore $c\neq 0$ and $(3c)^3\neq (2+c^3)^3$, so $$(2+c^3)^3-(3c)^3 = (c^3-1)^2(c^3+8) \neq 0.$$ Setting $a=b=1=-q$, we see that the algebra $A$ defined above is  $S(-c,-c,2)$. This algebra is a degenerate Sklyanin if and only if $c^3 = -8$, so $S(-c,-c,2)$ is a three-dimensional Sklyanin algebra. Moreover, $S(-c,-c,2)$ is also of type EC; for if $(6c^2)^3=(8-2c^3)^3$, then $$(6c^2)^3-(8-2c^3)^3 = 8(c^3-1)(c^3+8)^2 = 0,$$ which contradicts the assumption that $S(1,1,c)$ is of type EC. The automorphism of the point scheme of $S(-c,-c,2)$ is translation by $[-c:-c:2]$ in the group law with identity element $[1:-1:0]$, hence, no closed points of $\G(S(-c,-c,2))$ lie on $\Delta(\P^2)$. We conclude that for every point $$[\a_1:\b_1,\g_1] \times [\a_2: \b_2: \g_2] \in \G(S(-c,-c,2)),$$ if $r = \a_1 x+\b_1y+\g_1z$ and  $s = \a_2 x+\b_2y+\g_2z$ are elements in $S(1,1,c)_1$, then $r$ and $s$ are linearly independent and $rs = -sr$. Conversely, a pair $(r,s)$ of such $r$ and $s$, up to scaling in each component independently, determines a unique closed point in $\G(S(-c,-c,2))$. 
\end{proof}

\begin{rmk}
 The condition in Theorem \ref{not a skew polynomial ring} (2): $S(1,1,c)$ is of type EC is not necessary for the existence of stated $r$ and $s$. For example, if $c^3=-8$, then $S(1,1,c)$ is not of type EC, and $A = S(-c,-c,2)$ is a degenerate Sklyanin algebra. By \cite{Smith}, provided ${\rm char}\ \k\neq 3$, $A\cong \k\la u, v, w\ra/\la u^2, v^2, w^2\ra$. One checks that $A$ is semi-standard, and $\G(A)$ is a nonempty, non-identity relation on $E=\mathcal Z(uvw)$. Thus $\G(A)$ is not contained in $\Delta(\P^2)$, and one can use the isomorphism in \cite{Smith} to find linearly independent $r, s \in S(1,1,c)_1$ such that $rs=-sr$.
\end{rmk}

\begin{prop}
\label{counting S(1,1,c)}
Assume ${\rm char}\ \k\neq 2,3$ and let $j\in \k$. The number of distinct isomorphism classes of algebras $S(1,1,c)$ of type EC whose point scheme $E$ has j-invariant $j(E)=j$ is:
\begin{itemize}
\item three, if $j\neq 0, 12^3$,
\item one, if $j=0$, 
\item two, if $j=12^3$.
\end{itemize}
\end{prop}

\begin{proof}
Observing that the automorphism $\s$ of the point scheme of $S(1,1,c)$ is translation by $[1:1:c]$, a point of order 2, the proof is identical to that of \cite[Proposition 3.10]{Mat}, with the obvious necessary change from type B to type A.
\end{proof}

We complete our characterization of when an algebra $S(a,b,c)$ of type EC contains a two-dimensional AS-regular subalgebra by showing that no three-dimensional Sklyanin algebras contain a Jordan plane generated in degree 1.

\begin{prop}
\label{not a Jordan plane}
Let $S(a,b,c)$ be a three-dimensional Sklyanin algebra such that $abc\neq 0$. Assume ${\rm char}\ \k\neq 2, 3$. If $r, s\in S(a,b,c)_1$ satisfy $rs=sr-s^2$, then $r$ and $s$ are linearly dependent.
\end{prop}

\begin{proof}
Suppose $r, s\in S(a,b,c)_1$ satisfy $rs=sr-s^2$, and $abc\neq 0$. Without loss of generality, we assume $a=1$ and write $r=\a_1x+\b_1y+\g_1z$ and $s=\a_2x+\b_2y+\g_2z$ for $\a_i, \b_i, \g_i\in \k$. 
Examining the coefficients of $x^2$, $y^2$, $z^2$ in the equation $rs=sr-s^2$ implies the following identities hold in $\k$:
\begin{align*}
\a_2^2-c\b_1\g_2+c\g_1\b_2-c\b_2\g_2&=0\qquad (E_1)\\
\b_2^2-c\g_1\a_2+c\a_1\g_2-c\a_2\g_2&=0\qquad (E_2)\\
\g_2^2-c\a_1\b_2+c\b_1\a_2-c\a_2\b_2&=0\qquad (E_3).
\end{align*}
The linear combination $\a_2(E_1)+\b_2(E_2)+\g_2(E_3)$ yields
\begin{equation}
\label{eqn3}
\a_2^3+\b_2^3+\g_2^3-3c\a_2\b_2\g_2=0.
\end{equation}
Examining the coefficients of $zy$, $xz$, $yx$ in the equation $rs=sr-s^2$ implies:
\begin{align*}
(1+b)\g_1\b_2-(1+b)\b_1\g_2+(1-b)\b_2\g_2&=0\qquad (E'_1)\\
(1+b)\a_1\g_2-(1+b)\g_1\a_2+(1-b)\a_2\g_2&=0\qquad (E'_2)\\
(1+b)\b_1\a_2-(1+b)\a_1\b_2+(1-b)\a_2\b_2&=0\qquad (E'_3).
\end{align*}
The linear combinations $c(E'_i)-(1+b)(E_i)$ for $i=1,2,3$ show:
$$(1+b)\a_2^2=2c\b_2\g_2,\qquad
(1+b)\b_2^2=2c\a_2\g_2,\qquad
(1+b)\g_2^2=2c\a_2\b_2.$$
If $b=-1$, then since $c\neq 0$, these equations imply that at least two of $\a_2, \b_2, \g_2$ are zero. It then follows from $(E_1)$, $(E_2)$, $(E_3)$ that all of $\a_2, \b_2, \g_2$ are zero, so $s=0$. Suppose that $b+1\neq 0$. Note that the last three displayed equations show that $$(1+b)\a_2^3=2c\a_2\b_2\g_2 = (1+b)\b_2^3 = (1+b)\g_2^3.$$ Thus $\a_2^3=\b_2^3=\g_2^3$ and equation (\ref{eqn3}) implies
$3(1-b)\a_2^3=0.$ If $\a_2=0$, then $s=0$. If $b=1$, equations $(E'_1)$, $(E'_2)$, $(E'_3)$ simplify and show that all minors of 
$$\begin{bmatrix} \a_1 & \b_1 & \g_1\\ \a_2 & \b_2 & \g_2\\ \end{bmatrix}$$
vanish, whence $r$ and $s$ are linearly dependent.
\end{proof}

\subsection{Certain twisted tensor products as AS-regular algebras of type A}
We introduce a one-parameter family of algebras, $P(a)$, $a\in \k$. When $a\neq 1$, a member of this family is both a graded twisted tensor product of $\k_{-1}[r,s]$ and $\k[t]$, and an AS-regular algebra. When $a\neq 0,1$, $P(a)$ is of type EC, subtype A. We show that every algebra in the family $S(1,1,c)$ of Sklyanin algebras is isomorphic to some $P(a)$, and hence is a graded twisted tensor product of $\k_{-1}[r,s]$ and $\k[t]$.

Assuming that ${\rm char}\ \k\neq 2$, for any $a \in \k$, we define
$$P(a)=\dfrac{\k\la r,s,t\ra}{\la rs+sr, tr+rt-t^2-ar^2, t(r-s)+(r-s)t+2s^2\ra}.$$

\begin{prop}
\label{TTP conditions for P(a,b)}
The algebra $P(a)$ is a graded twisted tensor product of $\k_{-1}[r,s]$ and $\k[t]$ if and only if $a\neq 1$. Moreover, when $a\neq 1$, $r^2, s^2,$ and $t^2$ are central.
\end{prop}

\begin{proof}
Order monomials in the free algebra by left-lexicographic order via $r<s<t$. Then in degree 2, the Gr\"obner basis for the defining ideal of $P(a)$ is $\{sr+rs, ts-tr+st-rt-2s^2, t^2-tr-rt+ar^2\}.$
Resolving the overlaps $tsr$, $t^3$, and $t^2s$ via the Diamond Lemma yields two additional cubic Gr\"obner basis elements:
$$trs+tr^2-str+rtr+2rs^2 \qquad (1-a)tr^2-(1-a)r^2t.$$
If $a=1$, we have $\dim_{\k} P(a)_3=11$. However, any graded twisted tensor product of $\k_{-1}[r,s]$ and $\k[t]$ has the same Hilbert series as $\k_{-1}[r,s]\tsr_{\k} \k[t]$, namely $(1-z)^{-3}$. Thus $a\neq 1$ is necessary for $P(a)$ to be a graded twisted tensor product, and we assume this condition holds for the remainder of the proof. It is now straightforward to check that the newly-introduced overlaps $t^2rs$, $t^2r^2$, and $trsr$ resolve, hence
$$sr+rs, \qquad ts-tr+st-rt-2s^2, \qquad t^2-tr-rt+ar^2, $$
$$trs+r^2t-str+rtr+2rs^2, \qquad tr^2-r^2t $$
constitutes a Gr\"obner basis for the defining ideal of $P(a)$, and 
$$\{r^is^j(tr)^kt^{\e}\ |\ i,j,k\ge 0, \e\in\{0,1\} \}$$
is a $\k$-basis for $P(a)$. An easy counting argument shows the Hilbert series of $P(a)$ is $(1-z)^{-3}$. It is also easy to check that $r^2, s^2,$ and $t^2$ are central.

To complete the proof, we show that $\{ r^is^jt^k\ |\ i,j,k\ge 0\}$ is a spanning set for $P(a)$. Let $U$ be the $\k$-linear subspace of $P(a)$ spanned by $\{ r^is^jt^k\ |\ i,j,k\ge 0\}$. Given an arbitrary monomial $m \in P(a)$ in the generators $r, s, t$, we need to show that $m \in U$. Since $r$ and $s$ skew commute and $t^2$ is central, it suffices to show that $tr^is^j \in U$ for all $i, j \geq 0$. Observe that the defining relations show that $tr$ and $ts$ are both in $U$. Write $i = 2q+p$, $j = 2q'+p'$ for some integers $p, p', q, q'$ with $p, p' \in \{0,1\}$. Using the fact that $r^2$ and $s^2$ are central, we have $$tr^is^j = tr^{2q}r^p s^{2q'} s^{p'} = r^{2q} s^{2q'} tr^ps^{p'}.$$ If $p = 0$, it is clear that the last expression is in $U$. If $p = 1$, then note that 
$ trs^{p'} = (ar^2-rt+t^2)s^{p'} \equiv rts^{p'} \equiv 0 \pmod{U}.$
\end{proof}

\begin{prop}
\label{AS- regular for P(a,b)}
The algebra $P(a)$ is AS-regular if and only if $a\neq 1$.
\end{prop}

\begin{proof}
Let $S = P(a)$ and let $F=\k\la r, s, t\ra$ be the free algebra. The defining quadratic relations of $P(a)$ are: $$f_1 = rs+sr, \quad f_2 = tr+rt-t^2-ar^2, \quad f_3 = t(r-s)+(r-s)t+2s^2.$$ 

We have $M [r \ s \ t]^T = [f_1 \ f_2 \ f_3]^T$, where $$M = \begin{bmatrix} s & r & 0 \\ t-ar & 0 & r-t \\ t & 2s-t & r-s\end{bmatrix}.$$ 
Let $E$ be the closed subscheme of $\P^2$ determined by $\det M$. Recall that $S$ is nondegenerate provided that the rank of $M$ is equal to $2$ at every closed point of $E$. Let us write $M(p)$ for the matrix $M$ evaluated at a closed point $p \in E$. If $[r:s:t] \in E$, $r \ne 0$ or $s \ne 0$, and $r\ne s$ or $r\ne t$, then it is clear that ${\rm rank \ } M(r:s:t) = 2$. If $r = s = 0$, the second and third rows of $M(0:0:1)$ show that ${\rm rank \ } M(0:0:1) = 2$. At the point $[1:1:1] \in E$, $$M (1:1:1)= \begin{bmatrix} 1 & 1 & 0 \\ 1-a & 0 & 0 \\ 1 & 1 & 0\end{bmatrix},$$ so the condition $a\neq 1$ is equivalent to ${\rm rank \ } M(1:1:1) = 2$. We conclude that $S$ is nondegenerate if and only if $a\neq 1$.

Assume that $a\neq 1$. We now show that  $S$ is a standard algebra.  Let
$$
P=\begin{bmatrix}
a(s-r) & t-2s+ar & s-t\\
t-2s+ar & -2t+4s-2r & t-2s+ r\\
s-t & t-2s+ r & s-r\\
\end{bmatrix}
$$
and define $r_1, r_2, r_3\in F$ via $P[r\ s\ t]^T=[r_1\ r_2\ r_3]^T$. Note that $P=P^T$ so the algebra $F'=F/\la r_1, r_2, r_3\ra$ is standard. We shall prove that $F'=P(a)$. Since $P$ is symmetric, it suffices to show that the entries of $[r\ s\ t]P$ span the defining quadratic relations of $P(a)$. 
One checks that
$$[r\ s\ t]P=\begin{bmatrix}f_1 & f_2 & f_3\end{bmatrix}\begin{bmatrix} a & -2 & 1\\
1 & -1 & 0\\
-1 & 2 & -1
\end{bmatrix}.$$
Moreover, the determinant of the last matrix is equal to $a-1$. Hence $[r\ s\ t]P$ spans the quadratic relations of $P(a)$. Therefore $S$ is standard. We conclude by \cite[Theorem 1]{ATVI} that $S$ is AS-regular of dimension 3.
\end{proof}

\begin{prop}
\label{point scheme of P(a,b)}
 Assume ${\rm char}\ \k\neq 2, 3$. Let $a\in \k$ such that $a\neq 1$.
\begin{itemize}
\item[(1)] The point scheme of $P(a)$ is an elliptic curve if and only if $a \neq 0$. 
\item[(2)] The automorphism of the point scheme of $P(a)$ has order 2.
\item[(3)] Every elliptic curve, up to isomorphism, is the point scheme of some $P(a)$. 
\end{itemize}
\end{prop}

\begin{proof}
The condition $a \neq 1$, by Proposition \ref{AS- regular for P(a,b)}, guarantees that $P(a)$ is AS-regular, and so the point scheme of $P(a)$ and its associated automorphism are defined.
The matrix $M$ for $P(a)$ is given in the proof of Proposition \ref{AS- regular for P(a,b)}. One then checks that the point scheme of $P(a)$ is given by $E_{a}= \mathcal Z(F)$, where
$$F = \det M = ar^3-ar^2s-2rs^2+2rst+2s^2t-rt^2-st^2.$$ 

We note that $\frac{\partial F}{\partial t}=2(r+s)(s-t)$, and it follows that the only points satisfying $F=\frac{\partial F}{\partial t}=0$ are  $[0:0:1]$, $[1:1:1]$ and $[x:y:y]$ where $ax^2=y^2$. One readily checks that none of these points are singular when $a(a-1)\neq 0$. Conversely, if $a = 0$, then $[1:0:0]\in E_a$ is singular.  Therefore (1) holds.

The associated automorphism $\s_{a}$ of $E_{a}$ interchanges $[0:0:1]$ and $[1:1:1]$, and for $r\neq s$, computing the cross product of the first and third rows of $M$  shows
$$\s_{a}([r:s:t])=[r:-s:(2s^2-rt-st)/(r-s)].$$
It follows from an easy computation that $\s_{a}$ has order 2. 

For (3), suppose that $a(a-1) \neq 0$ so that the point scheme of $P(a)$ is the elliptic curve given by $F = 0$. A Weierstrass equation of this curve is 
$$Y^2Z+2XYZ-4(a-1)YZ^2=X^3+(4-a)X^2Z-4(a-1)XZ^2.$$
Using the formulae in \cite[III.1]{Sil} one checks that the $j$-invariant is given by
$$j_{a}=\frac{16(a^2+14a+1)^3}{a(a-1)^4}.$$
Since $\k$ is algebraically closed, this expression realizes every value in $\k$. It is well known that the $j$-invariant parametrizes elliptic curves up to isomorphism, see \cite[Proposition III.1.4(b)]{Sil} for example, hence (3) follows. 
\end{proof}

\begin{rmk}
\label{sigma_a is translation}
Suppose that $a\in \k$ satisfies $a(a-1) \neq 0$. Let $E_a$ denote the point scheme of the algebra $P(a)$. Then $E_a$ is an elliptic curve by Proposition \ref{point scheme of P(a,b)} (1).  Let us define, for use below, points $O = [0:0:1]$, $P = [1:1:1]$ of $E_a$, and take $O$ as the identity element for the group law on $E_a$, for all $a$. It is then straightforward, albeit tedious, to check that the automorphism $\s_a$ is translation by the point $P$.
\end{rmk}

\begin{thm}
\label{classifying P(a)}
If $a(a-1)\neq 0$, then $P(a)\cong P(a')$ if and only if $a' \in \{a, 1/a\}$.
\end{thm}

\begin{proof}
Assume that $a \in \k$ and $a(a-1)\neq 0$. Then $P(a)$ is AS-regular of type EC, and $P(a)$ is a geometric algebra. Clearly $P(a)\ncong P(a')$ if $a'\in\{0,1\}$, so suppose $P(a)\cong P(a')$ for $a'(a'-1)\neq 0.$
By Theorem \ref{geometric algebra isomorphism}, it suffices to characterize projective equivalences $\psi:E_a\to E_{a'}$ such that $\psi\s_a = \s_{a'}\psi$.  Let $\psi$ be such an equivalence.

The point $O = [0:0:1]$ is a flex point for both $E_a$ and $E_{a'}$. It is well known that a projective equivalence preserves the set of flex points and that the set of flex points is equal to the set of $3$-torsion points. Hence $\psi(O)$ is a 3-torsion point for $E_{a'}$. Since translation by any 3-torsion point extends to an automorphism of $\P^2$,  see \cite[Lemma 5.3]{Mori} for example, the map $\psi' = \t_{-\psi(O)} \psi: E_a \to E_{a'}$, where $\t_{-\psi(O)}$ is translation on $E_{a'}$ by the point $-\psi(O)$, is a projective equivalence. Moreover, by Remark \ref{sigma_a is translation}, $\s_{a'}$ is a translation and one checks that $\psi'\s_a = \s_{a'} \psi'$. Therefore there is no loss of generality in assuming $\psi(O)=O$, and that $\psi$ is represented (up to scale) by an invertible matrix of the form
$$\begin{bmatrix}
a_{11} & a_{12}  & 0\\
a_{21} & a_{22} & 0\\
a_{31} &  a_{32} & 1
\end{bmatrix}.$$
Since $\s_a([0:0:1])=[1:1:1]$, we see $\psi([1:1:1]) = [1:1:1]$ and it follows that $$a_{11}+a_{12}=a_{21}+a_{22}=a_{31}+a_{32}+1.$$ 
For the remainder of the proof we consider the action of $\psi$ on $p\in E_a-\{O,P\}.$
It is easy to check that $[0:0:1]$ and $[1:1:1]$ are the only points on $E_a$ with $r=s$, and $[0:0:1]$ is the only point where $r+s=0$. Thus $p=[r:s:t]$ with $(r-s)(r+s)\neq 0$. The equation $\psi\s_a = \s_{a'}\psi$ implies that there exists a scalar $K_{p} \neq 0$ such that
\begin{align*}
& a_{11}r-a_{12}s=K_{p}(a_{11}r+a_{12}s) &\qquad &(E_1) \\
& a_{21}r-a_{22}s = -K_{p}(a_{21}r+a_{22}s) &\qquad &(E_2) \\
& a_{31}r-a_{32}s+(2s^2-rt-st)/(r-s)=K_{p}(2S^2-RT-ST)/(R-S) &\qquad &(E_3),
\end{align*}
where $R=a_{11}r+a_{12}s,$ $S=a_{21}r+a_{22}s$, and $T=a_{31}r+a_{32}s+t$. (The conditions $a_{11}+a_{12} = a_{21}+a_{22}$ and $\psi$ is invertible ensure that $R-S \neq 0$.) Since $a_{11}+a_{12}=a_{21}+a_{22}$, subtracting $(E_2)$ from $(E_1)$ yields
$$(a_{11}-a_{21})(r+s)=K_{p}[(a_{11}+a_{21})r+(a_{12}+a_{22})s].$$
Since $\psi$ is invertible and $r+s\neq 0$, the left side is nonzero, hence
$$K_{p}=\dfrac{(a_{11}-a_{21})(r+s)}{(a_{11}+a_{21})r+(a_{12}+a_{22})s}.$$
Now equation $(E_1)$ can be rewritten as $a_{11}a_{21}r^2-a_{12}a_{22}s^2 =0$. Since this homogeneous equation holds for all $p\in E_a-\{O,P\}$, we have $a_{11}a_{21}=a_{12}a_{22}=0$. Thus either (i) $a_{12}=a_{21}=0$, $a_{11}=a_{22}\neq 0$, and $K_{p}=1$, or (ii) $a_{11}=a_{22}=0$, $a_{12}=a_{21}\neq 0$, and $K_{p}=-1$. 

In case (i), since equation $(E_3)$ holds for all $p\in E_a-\{O,P\}$, we see that $a_{31}=0$ and $a_{32} = a_{11}-1$. Then $\psi([r:s:t]) = [a_{11}r: a_{11}s: (a_{11}-1)s+t]$. Using the assumption that $[r:s:t] \in E_a$, the condition $\psi([r:s:t])\in E_{a'}$ is equivalent to $$(a'a_{11}^2-a)r^3-(a'a_{11}^2-a)r^2s-(a_{11}^2-1)rs^2+(a_{11}^2-1)s^3 = 0.$$ 
This equation holds for all $p\in E_a-\{O,P\}$ if and only if $a = a'$ and $a_{11}^2 = 1$. 

Similarly, in case (ii), considering $(E_3)$ yields $a_{31}=a_{12}$ and $a_{32}=-1$. So $\psi([r:s:t]) = [a_{12}s: a_{12}r: a_{12}r-s+t]$. Using $[r:s:t] \in E_a$, the condition $\phi([r:s:t])\in E_{a'}$ is equivalent to $$(a_{12}^2-a)r^3-(a_{12}^2-a)r^2s - (a'a_{12}-1)rs^2+ (a'a_{12}^2-1)s^3 = 0.$$ This equation holds for all $p\in E_a-\{O,P\}$ if and only if $a_{21}^2=a$ and $a'a=1$.
The result follows.
\end{proof}

\begin{cor}
\label{counting P(a)}
For $j\in \k$, the number of distinct isomorphism classes of algebras $P(a)$ such that $a(a-1)\neq 0$ and the point scheme $E_a$ has j-invariant $j_a=j$ is:
\begin{itemize}
\item[(1)] three, if $j\neq 0, 12^3$,
\item[(2)] one, if $j=0$, and
\item[(3)] two, if $j=12^3$.
\end{itemize}
\end{cor}

\begin{proof}
Recall that for $a(a-1)\neq 0$, the j-invariant of the elliptic curve $E_a$ is given by $j_{a}=16(a^2+14a+1)^3/a(a-1)^4.$
Observe that $j_a=j_{1/a}$, so solutions to $j_a=j$ occur in mutually inverse pairs. By Theorem \ref{classifying P(a)}, each pair corresponds to a single isomorphism class in the family $P(a)$.
 
If $j = 0$, then the equation $j_a = 0$ has two distinct solutions, so (2) follows.

If $j\neq 0$, define a polynomial $$f_j(x) = (x^2+14x+1)^3 - (j/16)x(x-1)^4 \in \k[t].$$ Observe that $j_a = j$ if and only if $a$ is a zero of $f_j(x)$. It is readily checked that
$$(x^2+14x+1)f'_j(x)-6(x+7)f_j(x)=(j/16)(x-1)^3(x+1)(x^2-34x+1).$$
Since $j\neq 0$ and $f_j(1)\neq 0$, if $f_j(x)$ has a multiple root at $x=a$, then we have $(a+1)(a^2-34a+1)=0$. For such an $a$, $f_j(a)=0$ implies $j=12^3$, so (1) follows. For (3), observe that $j_a=12^3$ implies $(a+1)^2(a^2-34a+1)^2=0$.
\end{proof}

AS-regular algebras whose point schemes are elliptic curves have recently been classified up to isomorphism by Itaba and Matsuno in characteristic 0. In \cite[Theorem 4.9]{IM}, the authors provide a classification of geometric algebras whose point schemes are elliptic curves. However, not all algebras listed in \cite[Theorem 4.9]{IM} are AS-regular. In \cite[Theorem 3.13]{Mat} (see also \cite[Remark 3.9]{Mat}) Matsuno completes the classification. We draw attention in particular to the family of algebras 
 $$R(a,b,c) = \frac{\k \la x, y, z \ra}{\la axz+bzy+cyx, azx+byz+cxy,ay^2+bx^2+cz^2\ra}$$
 where $[a:b:c]\in \P^2$ is a point on $E_{\l}=\mathcal Z(x^3+y^3+z^3-3\l xyz)$, $\l^3\neq 1$ such that $abc\neq 0$. When the $j$-invariant of $E_{\l}$ is not equal to $0$ or $12^3$, then these are the Type B algebras of \cite[Table 1]{Mat}; also see \cite[Theorem 4.9]{IM}.

\begin{lemma}\label{center of R(a,b,c)}
Let $[a:b:c]\in E_{\l}$, $\l^3\neq 1$, such that $abc \ne 0$. If $\dim_{\k} R(a,b,c)_3 = 10$, then $\dim_{\k} Z(R(a,b,c))_2 = 2$. 
\end{lemma}

\begin{proof}
We calculate a Gr\"obner basis to degree 3. Order monomials in the free algebra $\k\la x, y, z \ra$ using left-lexicographic order based on $x < y < z$. Without loss of generality assume that $c = 1$, and write the defining relations of $R(a,b,c)$ as: 
$$z^2+ay^2+bx^2\qquad zy+b^{-1}yx+ab^{-1}xz\qquad zx+a^{-1}byz+a^{-1}xy.$$
Then the overlaps in degree $3$ are: $z^3, z^2y, z^2x$. Resolving these yields: 
\begin{align*}
&a^{-2}(b^3-a^3)[y^2z-a^2b^{-2}x^2z+b^{-1}yxy-ab^{-2}xyx], \\
&b^{-2}(1-b^3)[yx^2-x^2y], \\
&a^{-2}(1-a^3)[xy^2-y^2x],
\end{align*}
respectively. 

There are twelve cubic monomials in the free algebra $\k\la x, y, z \ra$ that do not contain $z^2$, $zx$ or $zy$ as a submonomial. Assuming that $\dim_{\k} R(a,b,c)_3 = 10$, there must be exactly two linearly independent Gr\"obner basis elements in degree 3. Therefore exactly one of: $b^3-a^3$, $1-b^3$, $1-a^3$ must be zero. We consider cases.

\noindent{\bf Case 1:} $b^3-a^3 = 0\ne (1-b^3)(1-a^3)$. 

The Gr\"obner basis in degree 3 is spanned by: $yx^2-x^2y$ and $y^2x-xy^2$. It is straightforward to check that  $\{a^2b^{-2}x^2+y^2, xy+yx\}$ is a basis for $Z(R(a,b,c))_2$.

\noindent{\bf Case 2:} $1-b^3 = 0\ne (b^3-a^3)(1-a^3)$. 

The Gr\"obner basis in degree 3 is spanned by: $y^2z-a^2b^{-2}x^2z+b^{-1}yxy-ab^{-2}xyx$ and $y^2x-xy^2$. It is easy to check that $\{xy-axz+byz, y^2\}$ is a basis for $Z(R(a,b,c))_2$.

The remaining case is analogous to Case 2.
\end{proof}

\begin{thm}
\label{P(a,b) is an S(1,1,c)}
Let $\k$ be an algebraically closed field of characteristic $0$. When $a(a-1)\neq 0$, the twisted tensor product $P(a)$ is isomorphic as a graded algebra to a three-dimensional Sklyanin algebra $S(1,1,c)$ of type EC. Moreover, every algebra $S(1,1,c)$ of type EC is isomorphic to some $P(a)$ with $a(a-1)\neq 0$. 
\end{thm}

\begin{proof}
Assume $a(a-1)\neq 0$. Proposition \ref{TTP conditions for P(a,b)} ensures that $P(a)$ is a twisted tensor product of $\k_{-1}[x,y]$ and $\k[t]$. Moreover, by Proposition \ref{AS- regular for P(a,b)} and Proposition \ref{point scheme of P(a,b)}, $P(a)$ is AS-regular, its point scheme is an elliptic curve and the associated automorphism has order 2. 

Combining \cite[Theorem 4.9]{IM}, \cite[Remark 3.9]{Mat}, and \cite[Theorem 3.13]{Mat}, it follows that $P(a)$ is isomorphic to a Sklyanin algebra $S(a',b',c')$ or to one of the algebras $R(a',b',c')$ defined above. We note that the Type E and Type H algebras in \cite[Table 1]{Mat} are ruled out, as the automorphisms for those types do not have order 2. 

Since $\dim_{\k} P(a)_3=10$ and $\dim_{\k} Z(P(a))_2=3$, Lemma \ref{center of R(a,b,c)} implies $P(a)$ must be isomorphic to a three-dimensional Sklyanin algebra $S(a',b',c')$. Since $P(a)$ is of type EC, so is $S(a',b',c')$. Since $P(a)$ contains a subalgebra isomorphic to $\k_{-1}[r,s]$, Theorem \ref{not a skew polynomial ring}(1) implies $P(a)$ is isomorphic to $S(1,1,c)$.

Finally, the statement that every $S(1,1,c)$ of type EC is isomorphic to some $P(a)$ with $a(a-1) \neq 0$ follows directly from Proposition \ref{counting S(1,1,c)} and Corollary \ref{counting P(a)}.

\end{proof}

\bibliographystyle{plain}
\bibliography{bibliog2}

\end{document}